\documentclass[11pt,a4paper]{amsart}
\usepackage{amsmath,amssymb,amsthm,yhmath}
\usepackage[all]{xy}
\newcommand{\bo}{\boldsymbol}
\DeclareMathOperator{\CC}{CC}
\DeclareMathOperator{\cosec}{cosec}
\DeclareMathOperator{\daa}{\,/\!\!\!\!\mathsf A}
\DeclareMathOperator{\dual}{dual}
\DeclareMathOperator{\GL}{GL}
\DeclareMathOperator{\Gr}{Gr}

\DeclareMathOperator{\Id}{Id}
\DeclareMathOperator{\notI}{\,\backslash\!\!\mathsf I}
\DeclareMathOperator{\PV}{PV}
\DeclareMathOperator{\sgn}{sgn}
\DeclareMathOperator*{\wlim}{w-lim}
\newcommand{\bix}{\overline{X}}
\newcommand{\emptis}{\textsf{{\O}}}
\newcommand{\fli}{\mathsf T}
\newcommand{\opp}{\mathrm{opp}}
\newcommand{\zz}{\mathbb Z}
\newcommand{\zzh}{\,\backslash\!\!\!\!\mathbb Z}
\newcommand{\zzi}{\zz_i}
\newcommand{\zzj}{\zz_j}
\newcommand{\zzk}{\zz_k}
\newcommand{\zzl}{\zz_l}
\renewcommand{\[}{\begin{equation}\notag}
\renewcommand{\]}{\end{equation}}
\swapnumbers
\theoremstyle{definition}
\newtheorem{point}{}[section]
\newtheorem{conven}[point]{Convention}
\theoremstyle{plain}
\newtheorem{lemma}[point]{Lemma}
\newtheorem{theorem}[point]{Theorem}
\newcommand{\marginextend}[1]{ \addtolength{\oddsidemargin}{-#1}  \addtolength{\evensidemargin}{-#1}
  \addtolength{\textwidth}{#1}\addtolength{\textwidth}{#1}}
\newcommand{\updownextend}[1]{ \addtolength{\topmargin}{-#1}  \addtolength{\textheight}{#1}
\addtolength{\textheight}{#1}}
\marginextend{0.7cm}
\updownextend{0cm}
\author{Gyula Lakos}
\title{Linearization of skew-periodic loops and $\mathbb S^1$-cocycles}
\address{Alfr\'ed R\'enyi Institute of Mathematics, Hungarian Academy of Sciences,
Budapest, Pf. 127, H--1364, Hungary}
\address{Institute of Mathematics, E\"otv\"os University, P\'azm\'any P\'eter s.~1/C,
Budapest, H--1117, Hungary}
\email{lakos@cs.elte.hu}
\begin{document}
\maketitle
\begin{abstract}
We discuss linearization of skew-periodic loops.
We generalize the situation to linearization of non-commutative loops and $\mathbb S^1$-cocycles.
\end{abstract}
\section*{Introduction}
Linearization is a phenomenon when a nonlinear system is made linear at the cost
of adding  extra space-dimensions but retaining essential equivalence to the old system.
The best known example is that there is a homotopy
\[\begin{bmatrix}
a_0+a_1z+\ldots+a_nz^n&&&&&\\
&1&&&&\\
&&1&&&\\
&&&\ddots&&\\
&&&&1&\\
&&&&&1
\end{bmatrix}
\sim
\begin{bmatrix}
a_0&-z&&&&\\
a_1&1&-z&&&\\
a_2&&1&-z&&\\
\vdots&&&\ddots&\ddots&\\
a_{n-1}&&&&1&-z\\
a_n&&&&&1
\end{bmatrix}
\]
such that this homotopy is achieved via multiplication by invertible matrices.
In particular, if $a_0+a_1z+\ldots+a_nz^n$ was invertible then
the homotopy will be through invertible matrices.
In case of non-polynomial systems one has to use infinitely many dimensions
and it may not be clear how to proceed even for convergence reasons.

In terms of loops the linearization problem is like to transform an invertible
Laurent loop $a(z)=\sum_{n=-\infty}^\infty a_nz^n$ into a linear loop
$\vec a(z)=\vec a_0+\vec a_1 z$ by a homotopy.
In \cite{SS} a sufficiently \textit{nice} method was presented for the case just mentioned,
which is essentially the case of the pointed complex loops.

In this paper we look a bit further, we show how to discuss the case of real loops, and some other ones.
As it turns out the most general case is not necessarily the most enlightening, so
we will start with specific examples, and then we proceed further to the more abstract setting.
The most basic case is the case of what we call based loops,
ie. maps from $\mathbb S^1$ into the unit group of an algebra.
Having discussed this  case we turn to various generalizations.
One possible generalization is when we use non-commutative loops in the sense such that
the loop variable $\mathsf z$ would not commute with the coefficients $a_n$.
Another generalization is about linearization of $\mathbb S^1$-cocycles,
as explained below, which allows to treat generalized loops  which are like
maps from a principal $\mathbb S^1$-bundle into a unit algebra, or a regular $\mathbb S^1$-action.

Another thing is that we concentrate on linearizing skew-periodic loops.
Compared to the previous example this is like to try to linearize the  loop
$a(z)=\sum_{n\in\zz+\frac12} a_nz^n$ into the linear loop
$\vec a(z)=\vec a_{-1/2}z^{-1/2}+\vec a_{1/2} z^{1/2}$.
This is a more natural thing to do.
Our methods here are the essentially the same ones as in \cite{SS} but
we feel worthwhile to present the modifications necessary.

Let us summarize here how we expect to execute linearizatizations:
It is assumed to be done in two steps:
1. We add extra space-variables.
This is supposed to be done by a trivial extension.
2. We linearize by a homotopy.
Now, there are some requirements:
a.) If a loop was already linear then linearization should essentially do nothing.
b.) The construction should be sufficiently natural.

The case of based loops is somewhat misleading.
Having a base point is like having an anchor which allows linearization
of loops as such. But even there, what we really do is
the linearization of the cocycle $\tilde a(z,w)=a(z)a(w)^{-1}$.
In the more general case we should rather linearize cocycles obtained not necessarily from loops.
This is in accordance to the geometric picture: when we apply the clutching construction
for vector bundles then what we need is not necessarily a clutching function but a clutching cocycle
in terms of the structure group.
There is an other widely used way to represent cocycles: that is to take the
corresponding element from the skew-loop Grassmannian space (cf.~\cite{PS}).
As it turns out, the corresponding Grassmannian element directly refers to
the linearized cocycle. In particular, we do not loose information in the linearization process.

In Section 1 we consider some very basic matrix constructions associated to based loops.
In order to represent loops we consider skew-complex and real Fourier bases.
In Section 2 we define some auxiliary constructions, they are basically about
how to write down some matrices naturally.
In Section 3 we present some constructions which are more to the point:
some deformation matrices are introduced. Using them,
in Section 4  some very natural deformations of multiplication actions
are introduced.
Section 5 is a short digression about deformations of rapidly decreasing
matrices in similar manner.
In Section 6, at last, we present the linearizing homotopies for based loops.
By that time we will have more than enough information to see why these linearizations
work.
Section 7 is about generalization to non-commutative loops.
Section 8 is a more detailed discussion of transformation kernels on $\mathbb S^1$
which are useful alternatives to matrices, and which may motivate some of our constructions
in the next, last section:
Section 9 is about linearization of $\mathbb S^1$-cocycles.
It sort of embodies the naturality of the basic linearization construction.
As it yields a somewhat abstract setting it is useful to think through its possible
geometric interpretations.

Regarding analysis: We will
treat the smooth case, even if we make efforts to point into the direction of algebraic generalizations.
Here ``smooth'' is used in the sense of infinitely many differentiable or
rapidly decreasing (in case we consider a Fourier or matrix expansions).
People with geometric applications on their mind may know what modifications
are  required in specific cases.

Some of the length of this paper is gained from that
we try to show the inner workings of the constructions.
In particular, most of our constructions come in three varieties:
the left perturbation case $(L)$, the ``unitary'' case $(\emptis)$, and
the right perturbation case $(R)$.
(The polynomial linearizing homotopy above would be of type ``$L$''.)
It would not be hard to give a much shorter exposition of the ``unitary'' case
alone, but some generalizations fit more to the ``$L$'' or ``$R$'' cases.

In terms of proofs, some of them are tedious matrix computations.
The author has decided not to include them here but to
present those key points which give the frame of the computations making them straightforward.

The author would like to thank for the hospitality of the Alfr\'ed R\'enyi Institute of Mathematics, and
in particular to P\'eter P\'al P\'alffy.

\section{Matrices associated to based loops}
\begin{point}
In what follows $\mathfrak A$ will be a locally convex algebra, $1\in\mathfrak A$,
although it is not hard to transcribe the constructions to the non-unital case.
We can define the opposite algebra $\mathfrak A^{\opp}$ such that
to any  $X\in\mathfrak A$ we  associate the opposite element  $X^\opp\in\mathfrak A^\opp$.
This operation is linear but opposite elements multiply in a contravariant manner:
$X^\opp\cdot Y^\opp=(YX)^\opp$.

A more usual manifestation of this operation is when
$\mathfrak A$ contains a distinguished skew-involution $g$ and there is
an isomorphism $\psi:\mathfrak A^\opp\rightarrow\mathfrak A$ given such that $\psi(g^\opp)=-g$.
Then the star operation  $a^*=\psi(a^\opp)$ can be used, it yields a $g$-antilinear operation.
\end{point}

\begin{point}\textit{Notation.}
In order to simplify the notation,  we make the following convention:
In the associated constructions instead of ``$x\mapsto a(x)$'' we will often write simply ``$a$''.
This assumes that the variable of $a$ is $x$,
and we do this only if there is no danger of confusion.

Matrices will be considered as linear combinations of elementary matrices.
General elementary matrices $\mathbf e_{n,m}=\mathbf e_n\otimes\mathbf  e_m^\top$ are
considered as tensor products of elementary column  and row  matrices.
In what follows ``$\ddots$'' denotes a diagonal pattern of something
and ``$\pmb\ddots$'' denotes a diagonal pattern of $1$'s.
If $S$ is an appropriate set of indices then we define the identity, flip, and shift matrices as
\[\mathbf 1_{S}=\sum_{s\in S}\mathbf e_{s,s},\qquad \fli_{-S,S}=\sum_{s\in S}\mathbf e_{-s,s},
\qquad \mathsf Z^n_{n+S,S}=\sum_{k\in S}\mathbf e_{n+k,k},\]
respectively.
If two lower indices are the same then only one is used.
We define $S^+=\{s\in S\,:\,s>0\}$,  $S^{++}=\{s\in S\,:\,s>1\}$, etc.
We will mainly be interested in the cases $S=\zz$ and
$S=\zzh\equiv \zz+\frac12$.
In case we do not fix which one we use we apply the notation $\zzi$,
$i\in\bo Z_2=\mathbb Z/2\mathbb Z$ where $\zz_0=\zz$ and  $\zz_1=\zzh$.
The complementer set is $\zzi'=(\frac12\zz)\setminus \zzi=\zzi+\frac12$.
When writing  $\zz\times$ matrices  we draw two lines, above and below the $0$th row,
when writing  $\zzh\times$ matrices  we draw one line between the $\tfrac12$th
and $(-\tfrac12)$th rows; similarly with columns.

If $\phi:\mathfrak A\rightarrow\mathfrak A'$ is a homomorphism and $\mathbf A$ is a
$\mathfrak A$-valued matrix then let simply  $\phi(\mathbf A)$ denote the matrix which
is obtained from $\mathbf A$ by replacing every entry $a_{n,m}$ with its homomorphic image $\phi(a_{n,m})$.
The generalization of taking opposites is given by taking $\mathbf A^{\top\opp}$,
ie. taking the transpose of $\mathbf A$ and replacing every entry $a_{m,n}$ by $a_{m,n}^\opp$.
\end{point}

\begin{point}\textit{Loops and Fourier expansions.}
Let $\mathfrak B^0=\mathcal C_{\zz}(\mathfrak A)$ be the space of smooth periodic $\mathfrak A$-valued loops,
ie.~of functions $a:\mathbb R\rightarrow\mathfrak A$ such that $a(x+2\pi)=a(x)$.
Similarly, let $\mathfrak B^1=\mathcal C_{\zzh}(\mathfrak A)$ be the space of smooth
skew-periodic $\mathfrak A$-valued loops,
ie.~of functions $a:\mathbb R\rightarrow\mathfrak A$  such that $a(x+2\pi)=-a(x)$.
Then $\mathfrak B=\mathfrak B^0\oplus \mathfrak B^1$ is a naturally $\bo Z_2$-graded algebra.
Such loops allow Fourier expansions in various ways:

a.) Suppose that $g$ is a skew-involution. Let us define $z(x)=\cos x+g\sin x$ where $x\in\mathbb R$.
We use the notation $z^\alpha(x)=\cos \alpha x+g\sin \alpha x$.
Then $a\in\mathfrak B^i$ allows a Fourier expansion $a(x)=\sum_{n\in\zzi}z^n(x)  a_n$.
Using the correspondence $\mathbf e_{n}\leftrightarrow z^n(x)$ of bases we can define the column matrix
\[\mathsf S_{\zzi}^{g}(a)=\sum_{n\in\zzi}a_n\mathbf e_n.\]

Similarly, we  can also define the row matrices $\bar{\mathsf S}_{\zzi}^{g}(a)$:
For a skew-involution $g$ we will use the short notation
$[K]^+_g=\frac{K-gKg}2$,  $[K]^-_g=\frac{K+gKg}2$,
for the $g$-invariant and   $g$-skew-invariant parts of $K$, respectively.
We see that $[K]^+_gz^\alpha(x)=z^\alpha(x)[K]^+_g$ but $[K]^-_gz^\alpha(x)=z^{-\alpha}(x)[K]^-_g$.
So, another form of the Fourier expansion is
$a(x)=\sum_{n\in\zzi}\left( [a_n]^+_g +[a_{-n}]^-_g \right) z^n(x)$.
Then, we let
\[\bar{\mathsf S}_{\zzi}^{g}(a)=\sum_{n\in\zzi}\left([a_n]^+_g+[a_{-n}]^-_g\right)\mathbf e_{-n}^\top.\]

For $K\in\mathfrak A$ it yields $\mathsf S_{\zzi}^{g}(aK)=\mathsf S_{\zzi}^{g}(a)K$
and $\mathsf S_{\zzi}^{g}(Ka)=([K]^+_g\mathbf 1_{\zzi}+[K]^-_g\fli_{\zzi} )\mathsf S_{\zzi}^{g}(a)$,
ie.~it is right $\mathfrak A$-linear and left  $\mathfrak A$-quasi-linear.
Similarly, $\bar{\mathsf S}_{\zzi}^{g}(Ka)=K\bar{\mathsf S}_{\zzi}^{g}(a)$
and $\mathsf S_{\zzi}^{g}(aK)=\bar{\mathsf S}_{\zzi}^{g}(a)([K]^+_g\mathbf 1_{\zzi}+[K]^-_g\fli_{\zzi} )$,
ie.~it is left $\mathfrak A$-linear and right $\mathfrak A$-quasi-linear.

When appropriate, right $\mathfrak A$-invariant linear maps from $\mathfrak B^j$ to $\mathfrak B^i$
can be represented by $\zzi\times\zzj$-matrices acting of the
left of the column vectors considered above. If the action is $A$ then for the representing matrix
$\mathbf A^g$ we have $ \mathsf S_{\zzj}^{g}(Aa)=\mathbf A^g\mathsf S_{\zzi}^{g}(a)$.
Similarly, through the representation $\bar{\mathsf S}_{\zzi}^{g}(a)$
such matrices can also be used to represent left $\mathfrak A$-linear actions
form $\mathfrak B^i$ to $\mathfrak B^j$.

b.) The problem with the skew-complex Fourier bases above is that sometimes
there are no  skew-involutions $g$ in $\mathfrak A$ at all.
For arithmetic purposes it is always possible to  adjoin an imaginary element $\mathrm i$
to $\mathfrak A$, but this might change the invariants associated to $\mathfrak A$.
Hence it is much more natural to use real Fourier series.
In the real representations we use bases corresponding to
$\mathbf e_n\leftrightarrow \sqrt2 \cos nx$ ($n>0$),
$\mathbf e_0\leftrightarrow 1$,
$\mathbf e_{-n}\leftrightarrow \sqrt2 \sin nx$ ($n>0$).
Our formulas using expansions in terms of  $z(x)=\cos x +g\sin x$ can be transcribed to the real base
with  the help of the  matrices
\[\mathsf N^{\mathrm r,g}_{\zz}=\mathbf e_{0,0}+\frac{1}{\sqrt2}\sum_{n\in \zz^+}
-g\mathbf e_{-n,-n}+g\mathbf e_{-n,n}+\mathbf e_{n,-n}+\mathbf e_{n,n},\]
\[\mathsf N^{\mathrm r,g}_{\zzh}=\frac{1}{\sqrt2}\sum_{n\in \zzh^+}
-g\mathbf e_{-n,-n}+g\mathbf e_{-n,n}+\mathbf e_{n,-n}+\mathbf e_{n,n},\]
and
$\mathsf N^{g,\mathrm r}_{\zzi}=(\mathsf N^{\mathrm r,g}_{\zzi})^{-1}
=(\mathsf N^{\mathrm r,-g}_{\zzi})^\top$.
Indeed, using the superscript ``$\mathrm r$'' in order to indicate real expansion we find that
$\mathsf S_{\zzi}^{\mathrm r}(a)
=\mathsf N^{\mathrm r,g}_{\zzi}\mathsf S_{\zzi}^{g}(a)$.
Similarly for $\bar{\mathsf S}_{\zzi}^{\mathrm r}(a)$.
These real representations are plainly left and right $\mathfrak A$-linear.
Furthermore, if $\mathbf A^{g}$ was a $\zzi\times\zzj$
matrix representing an action then the corresponding matrix in the real basis yields
\[\mathbf A^{\mathrm r}=\mathsf N^{\mathrm r,g}_{\zzi}\mathbf A^{g}(\mathsf N^{\mathrm r,g}_{\zzj})^{-1}.\]

In practice, it is easier to carry out computations in the skew-complex representation,
especially with  $g=\mathrm i$, then translate the result to the real case.
If the computation was not clear for $g\neq\mathrm i$ then we can still obtain the
corresponding formula via the real case.
\end{point}

\begin{point}\textit{Multiplication actions.}
Suppose that the  smooth function $x\mapsto a(x)$ is in $\mathfrak B^l$ and
it has Fourier expansion ${a(x)=\sum_{n\in\zzl}z^n(x)  a_n}$.
Suppose that $\zzi=\zzl+\zzj$, then $a\mathfrak B^j\subset\mathfrak B^i$.
Now, as $a(x) z^k(x)=\sum_{n\in\zzl}z^{n+k}(x)[a_n]^+_g+z^{n-k}(x)[a_n]^-_g$
we see that multiplication  by $a(x)$ on the left can be represented by the matrix
\[\mathsf U_{\zzi,\zzj}^{g}(a)=
\sum_{n\in\zzl,k\in\zzj}[a_n]^+_g\mathbf e_{n+k,k}+[a_n]^-_g\mathbf e_{n-k,k}=
\sum_{n\in\zzl}[a_n]^+_g\mathsf Z^n_{\zzi,\zzj}
+[a_n]^-_g\mathsf Z^n_{\zzi,\zzj}\fli_{\zzj}.\]
We see that $\mathsf S_{\zzi}^g(a)=\mathsf U_{\zzi,\zz}^g(a)\mathbf e_0$,
$\bar{\mathsf S}_{\zzj}^g(a)=\mathbf e_0^\top\mathsf U_{\zz,\zzj}^g(a)$.
These skew-complex multiplication representations are $\mathfrak A$-quasi-linear
in the sense that for $K\in\mathfrak A$
\[\mathsf U_{\zzi,\zzj}^{g}(Ka)=([K]^+_g\mathbf 1_{\zzi}+[K]^-_g\fli_{\zzi} )\mathsf U_{\zzi,\zzj}^{g}(a),\]
\[\mathsf U_{\zzi,\zzj}^{g}(aK)=\mathsf U_{\zzi,\zzj}^{g}(a)([K]^+_g\mathbf 1_{\zzj}+[K]^-_g\fli_{\zzj} ).\]

Furthermore, the representations above are multiplicative: ie.~the identity
\[\mathsf U^{g}_{\zzi,\zzj}(a_1)\mathsf U^{g}_{\zzj,\zzk}(a_2)=\mathsf U^{g}_{\zzi,\zzk}(a_1a_2)\]
holds.
Another simple but very useful observation which we call ``real multiplication commutativity''
is that for a periodic real scalar function $b$ we know that
\[\mathsf U^g_{\zzi,\zzj}(a)\mathsf U^g_{\zzj}(b)=\mathsf U^g_{\zzi}(b)\mathsf U^g_{\zzi,\zzj}(a).\]

Everything can be transcribed to the real base, yielding various matrices
$\mathsf U^{\mathrm r}_{\zzi,\zzj}(a)$. These representations are plainly
$\mathfrak A$-linear.
\end{point}
\begin{conven}
If we use indices $i,j,l$ then it with be understood that $\zzi=\zzl+\zzj$.
\end{conven}
\begin{point}
The orientation reversal action $C$ is defined by $(x\mapsto a(x)) \mapsto (x\mapsto a(-x))$.
If $a\in \mathfrak B^i$ and $a(x)=\sum_{n\in\zzi} z^n(x) a_n $ then $Ca(x)=\sum_{n\in\zzi} z^{-n}(x)a_n $.
Its matrix is $\mathsf C^g_{\zzi}=\fli_{\zzi}$ in the skew-complex base, and
$\mathsf C^{\mathrm r}_{\zzi}=-\mathbf 1_{\zzi^-}+\mathbf 1_{\zzi\setminus\zzi^-}$ in the real base.
\end{point}

\begin{conven}
Instead of a pair of indices ``$g,x_1$'' we may write ``$z_1$'' where $z_1=\cos x_1+g\sin x_1$.
We do this only when there is no danger of confusion.
\end{conven}

\begin{point} \textit{$\delta$-functions.}
Let $x_1$ be a constant or a variable independent from $x$.
They are not smooth functions  but only distributions, nevertheless
we can consider the periodic and skew-periodic ``algebraic'' delta functions
\[\delta_{\zz}^{x_1}(x)=2\pi\sum_{n\in\zz} \delta_{x_1+2\pi n}(x)\qquad\text{and}\qquad
\delta_{\zzh}^{x_1}(x)=2\pi\sum_{n\in\zz} (-1)^n\delta_{x_1+2\pi n}(x).\]

The column and row matrices representing  these distributions are
\[\boldsymbol\delta^{g,x_1}_{\zzi}
\equiv\boldsymbol\delta^{z_1}_{\zzi}
=\sum_{n\in\zzi}z_1^{-n}\mathbf e_{n},\qquad
\bar{\boldsymbol\delta}^{g,x_1}_{\zzi}
\equiv\bar{\boldsymbol\delta}^{z_1}_{\zzi}
=\sum_{n\in\zzi}z_1^{n}\mathbf e_{n}^\top.\]
The real-based matrices are
\[\boldsymbol\delta^{\mathrm r,x_1}_{\zz}=
\sum_{n\in\zz^+}\mathbf e_{-n}\sqrt2\sin nx_1
+\mathbf e_0+\sum_{n\in\zz^+}\mathbf e_{n}\sqrt2\cos nx_1,\qquad
\bar{\boldsymbol\delta}^{\mathrm r,x_1}_{\zz}
=(\boldsymbol\delta^{\mathrm r,x_1}_{\zz})^\top,\]
\[\boldsymbol\delta^{\mathrm r,x_1}_{\zzh}=
\sum_{n\in\zzh^+}\mathbf e_{-n}\sqrt2\sin nx_1
+\sum_{n\in\zzh^+}\mathbf e_{n}\sqrt2\cos nx_1,\qquad
\bar{\boldsymbol\delta}^{\mathrm r,x_1}_{\zzh}=
(\boldsymbol\delta^{\mathrm r,x_1}_{\zzh})^\top.\]
\end{point}

\begin{point}\textit{Rotation actions.}
Suppose that $x_1$ is a constant or a variable independent from $x$.
The rotation action $R^{x_1}$ will be given as $(x\mapsto a(x)) \mapsto (x\mapsto a(x-x_1))$ for us.
We set $z_1=\cos x_1+g\sin x_1$. We define the rotation matrices as
\[\mathsf R^{g,x_1}_{\zzi}\equiv \mathsf R^{z_1}_{\zzi}=\sum_{n\in\zzi}z_1^{-n}\mathbf e_{n,n}.\]
Indeed,
$\mathsf R^{g,x_1}_{\zzi}\mathsf S_{\zzi}^{g}(a)=\mathsf S_{\zzi}^{g}(R^{x_1}a)$, and
$\bar{\mathsf S}_{\zzi}^{g}(a)(\mathsf R^{g,x_1}_{\zzi})^{-1}=\bar{\mathsf S}_{\zzi}^{g}(R^{x_1}a)$.
More generally,
\[\mathsf R^{g,x_1}_{\zzi}\mathsf U_{\zzi,\zzj}^{g}(a)(\mathsf R^{g,x_1}_{\zzj})^{-1}
=\mathsf U_{\zzi,\zzj}^{g}(R^{x_1}a).\]
Also, $\mathsf R^{g,x_1}_{\zzi}\boldsymbol\delta^{g,x_2}_{\zzi}
=\boldsymbol\delta^{g,x_1+x_2}_{\zzi}$,
$\bar{\boldsymbol\delta}^{g,x_1}_{\zzi}(\mathsf R^{g,x_2}_{\zzi})^{-1}
=\bar{\boldsymbol\delta}^{g,x_1+x_2}_{\zzi}$.
The real-based matrices are
\[\mathsf R_{\zz}^{\mathrm r,x_1}=\mathbf e_{0,0}+\sum_{n\in \zz^+}
\mathbf e_{-n,-n}\cos nx_1+\mathbf e_{-n,n}\sin nx_1-\mathbf e_{n,-n}\sin nx_1+\mathbf e_{n,n}\cos nx_1,\]
\[\mathsf R_{\zzh}^{\mathrm r,x_1}=\sum_{n\in \zzh^+}
\mathbf e_{-n,-n}\cos nx_1+\mathbf e_{-n,n}\sin nx_1-\mathbf e_{n,-n}\sin nx_1+\mathbf e_{n,n}\cos nx_1.\]
\end{point}

\begin{point} \textit{Principal values of simple poles.}
Also quite important distributions are the distributions of Cauchy principal values of simple
poles. The periodic versions are
\[h_{\zz}^{x_1}(x)=\PV\cot\frac{x-x_1}2 \qquad\text{and}\qquad
h_{\zzh}^{x_1}(x)=\PV\cosec\frac{x-x_1}2 .\]

The column and row matrices representing these distributions are
\[\mathsf h^{g,x_1}_{\zzi}
\equiv\mathsf h^{z_1}_{\zzi}
=\sum_{n\in\zzi}-g(\sgn n)z_1^{-n}\mathbf e_{n},\qquad
\bar{\mathsf h}^{g,x_1}_{\zzi}
\equiv\bar{\mathsf h}^{z_1}_{\zzi}
=\sum_{n\in\zzi}-g(\sgn n)z_1^{n}\mathbf e_{n}^\top.\]
The real-based matrices are:
\[\mathsf h^{\mathrm r,x_1}_{\zzi}
=\sum_{n\in\zzi^+}\mathbf e_{-n}\sqrt2\cos nx_1 -\sum_{n\in\zzi^+}\mathbf e_n\sqrt2\sin nx_1,\qquad
\bar{\mathsf h}^{\mathrm r,x_1}_{\zzi}=-(\mathsf h^{\mathrm r,x_1}_{\zzi})^\top.\]
\end{point}

\begin{point}\textit{Integration and Hilbert transforms.}
For $u\in\mathfrak A$ the matrices
\[\textstyle{\int}^g_{\zz}(u)=\textstyle{\int}^{\mathrm r}_{\zz}(u)=u\mathbf e_{0,0}\]
represent the ``linear integration'' operation $\int_{\zz}: a\mapsto u\int_{y=0}^{2\pi} a(y)\frac{dy}{2\pi}$
on $\mathfrak B^0$.
Let
\[\mathsf H^g_{\zzi}=\sum_{n\in\zzi}-g (\sgn n)\mathbf e_{n,n}
=g\mathbf 1_{\zzi^-}-g\mathbf 1_{\zzi^+},
\]\[
\mathsf H^g_{\zz}[G]=g\mathbf 1_{\zz^-}+G\mathbf e_{0,0}-g\mathbf 1_{\zz^+}.\]

Then $\mathsf H^g_{\zz}$ and $\mathsf H^g_{\zzh}$
are the matrices of the Hilbert transforms $H_{\zz}$ and $H_{\zzh}$ on $\mathfrak B^0$ and $\mathfrak B^1$,
respectively. Note that $H_{\zz}$  is not invertible while $H_{\zzh}$ is invertible.
If $G$ is a skew-involution then a reasonable invertible
substitute for the even Hilbert transform is $\mathsf H^g_{\zz}[G]$.
The real-based matrices are
\[\mathsf H_{\zzi}^{\mathrm r}=\sum_{n\in\zzi}\sgn n \,\mathbf e_{-n,n}\qquad\text{and}\qquad
\mathsf H_{\zz}^{\mathrm r}[G]=G\mathbf e_{0,0}+\sum_{n\in\zz}\sgn n \,\mathbf e_{-n,n}.\]
\end{point}

\begin{point}\textit{Trivial extensions of loops.}
Let $z_1=\cos x_1+g\sin x_1$ again, and $u$ arbitrary.
We define the matrices
\[\mathsf L_{\zz}^{g,x_1}\equiv
\mathsf L_{\zz}^{z_1}=z_1\mathbf 1_{\zz^-}+\mathbf e_{0,0}+z_1^{-1}\mathbf 1_{\zz^+}
=\mathbf e_{0,0}+\mathbf 1_{\zz\setminus\{0\}}\cos x_1+\mathsf H_{\zz}^{g}\sin x_1,\]
\[\mathsf L_{\zzh}^{g,x_1}\equiv
\mathsf L_{\zzh}^{z_1}=z_1\mathbf 1_{\zzh^-}+z_1^{-1}\mathbf 1_{\zzh^+}
=\mathbf 1_{\zzh}\cos x_1+\mathsf H_{\zzh}^{g}\sin x_1,\]
\[\mathsf L_{\zz}^{g,x_1}[u]\equiv
\mathsf L_{\zz}^{z_1}[u]=z_1\mathbf 1_{\zz^-}+u\mathbf e_{0,0}+z_1^{-1}\mathbf 1_{\zz^+}
=u\mathbf e_{0,0}+\mathbf 1_{\zz\setminus\{0\}}\cos x_1+\mathsf H_{\zz}^{g}\sin x_1.\]
If $u$ is a general loop with variable $x_1$ then
$\mathsf L_{\zz}^{g,x_1}[u]\equiv\mathsf L_{\zz}^{z_1}[u]$ can be
considered as the trivial extension of $u$ by infinitely many space dimensions.
The real-based matrices can quickly be recovered from the various constituents.
\end{point}

\begin{point}
Up to this point we have defined several matrices $\mathsf M^g$.
We have seen that these matrices allow real-based transcriptions $\mathsf M^{\mathrm r}$
which do not depend on any skew-involution $g$, and, in fact, they do not require
the existence of any such skew-involution.
Colloquially, $\mathsf M^g$ allows a real version  $\mathsf M^{\mathrm r}$.
In what follows we use the superscript ``$g$'' in order to indicate that such a
transcription is possible, even if we do not elaborate it in detail.
Then, it goes without saying that any formula involving various matrices $\mathsf M^g$, linear combinations
and multiplications allows a similar formula involving corresponding matrices $\mathsf M^{\mathrm r}$.
\end{point}
\begin{point}\textit{Algebraic kernel functions.}
For a rapidly decreasing $\zzi\times\zzj$ matrix
\[\mathbf A^g=\sum_{n\in\zzi,m\in\zzj}a_{n,m}\mathbf e_{n,m}\]
we define the its algebraic kernel function as
\[A(x_1,x_2)=\bar{\bo\delta}^{g,x_1}_{\zzi}\mathbf A^g\bo\delta^{g,x_2}_{\zzi}
=\sum_{n\in\zzi,m\in\zzj}z_1^na_{n,m}z_2^{-m}.\]

If $\mathbf A^{\mathrm r}$ is the real-based transcription of $\mathbf A^g$ then
$A(x_1,x_2)=\bar{\bo\delta}^{\mathrm r,x_1}_{\zzi}\mathbf A^{\mathrm r}\bo\delta^{\mathrm r,x_2}_{\zzi}$.
In particular, if $A(x_1,x_2)$ is an expression independent of $g$ then $\mathbf A^g$ allows a
\textit{natural} real transcription.
In fact, instead of using (and examining) real matrices it is simpler to use the algebraic kernel function
$A(x_1,x_2)$.
Transformation kernels can also be defined, more generally, for matrices of temperate growth
but then they yield distributions.
\end{point}
\begin{point}\textit{Grassmannians.} The loop space skew-Grassmannians are defined by
\[\Gr^g_{\zzh} \mathfrak B^0=\{\mathsf U^{g}_{\zzh}(a) \mathsf H_{\zzh}^g
\mathsf U^{g}_{\zzh}(a^{-1})\,:\,a\in\mathfrak B^0 \text{ is invertible}\},\]
\[\Gr^g_{\zz} \mathfrak B^1=\{\mathsf U^{g}_{\zz,\zzh}(a) \mathsf H_{\zzh}^g
\mathsf U^{g}_{\zzh,\zz}(a^{-1})\,:\,a\in\mathfrak B^1 \text{ is invertible}\}.\]

We use the notations $\mathsf J^g_{\zzh}(a)$ and $\mathsf J^g_{\zz}(a)$
for the skew-involutions associated to the various invertible loops $a$ as above, respectively.
\end{point}
\begin{theorem}
$\mathsf J^g_{\zzh}$ can be ``inverted'' by the formula
\[a(x_1)a(x_2)^{-1}=\bar{\boldsymbol\delta}_{\zz}^{g,x_1}\biggl(
\mathsf U^{g}_{\zz,\zzh}(\sin\tfrac x2)\mathsf J_{\zzh}^{g}(a)\mathsf U^{g}_{\zz,\zzh}(\cos\tfrac x2)-
\mathsf U^{g}_{\zz,\zzh}(\cos\tfrac x2)\mathsf J_{\zzh}^{g}(a)\mathsf U^{g}_{\zz,\zzh}(\sin\tfrac x2)
\biggr){\boldsymbol\delta}_{\zz}^{g,x_2}.\notag\]

Similar formula holds with $\mathsf J^g_{\zz}$ but with $\zz$ and $\zzh$ interchanged.
\begin{proof} For $g=\mathrm i$ the inversion formula
\begin{multline}
a(x_1)a(x_2)^{-1}=
\left(\sum_{n\in\zz}  z_1^n\mathbf e_n^\top \right)
\mathsf U_{\zz}^g(a)\mathbf e_{0,0}\mathsf U_{\zz}^g(a^{-1})
\left(\sum_{m\in\zz}  z_2^{-m}\mathbf e_{m} \right)=\\
=\bar{\boldsymbol\delta}^{z_1}_{\zz}
\frac1{2\mathrm i}\biggl(\mathsf U^g_{\zz,\zzh}( z^{1/2})
\mathsf J_{\zzh}(a)\mathsf U^g_{\zzh,\zz}(z^{-1/2})
-\mathsf U^g_{\zz,\zzh}(z^{-1/2})
\mathsf J_{\zzh}(a)\mathsf U^g_{\zzh,\zz}( z^{1/2})\biggr)
\boldsymbol\delta^{z_2}_{\zz}
\notag\end{multline}
is easy to check.
Expanded further, it yields our result, which form remains valid even after changing bases.
\end{proof}
\end{theorem}
\begin{theorem} The following are equivalent:

i.) $J\in \Gr^g_{\zz} \mathfrak B^1$

ii.) $J$ is a rapidly decreasing perturbation of $\mathsf H_{\zz}^{g}$, $J^2=-\mathbf 1_{\zz}$, and
\[\mathsf U^{g}_{\zz}(\cos x)+J\mathsf U^{g}_{\zz}(\cos x) J=
\mathsf U^{g}_{\zz}(\sin x)J-J\mathsf U^{g}_{\zz}(\sin x).\]

Similar statement holds for $J\in \Gr^g_{\zzh} \mathfrak B^0$.
\begin{proof}
Assume that $g=\mathrm i$. Now, it is not hard to see that (i) holds
if and only if along the rapid decrease and square conditions also the conditions
\[\frac{1-\mathrm i J}{2}\mathsf U^g_{\zz}(z)\frac{1+\mathrm i J}{2}=0\qquad\text{and}\qquad
\frac{1+\mathrm i J}{2}\mathsf U^g_{\zz}(z^{-1})\frac{1-\mathrm i J}{2}=0\]
hold. (In an unprecise but possibly familiar language this is just the requirement that
the $(-\mathrm i)$-eigenspace of $J$ is $z$-invariant for multiplication, and the
$\mathrm i$-eigenspace of $J$ is $z^{-1}$-invariant for multiplication.)
Expanding the two conditions in terms of $z=\cos x+\mathrm i\sin x $
we obtain, somewhat surprisingly, a single equation as in (ii) for the case $g=\mathrm i$.
Changing bases the result follows in full generality.
\end{proof}
\end{theorem}
\begin{lemma} If $G$ is a skew-involution then
\[\mathsf H_{\zz}^{g}[G]=\mathsf J_{\zz}^{g}(\cos\tfrac x2+G\sin\tfrac x2).\]
\qed\end{lemma}

\begin{point}\textit{Poisson and Hilbert-Poisson transforms.}
Suppose that
$a=\sum_{n\in\zz} z^n a_n $. For $-1\leq r\leq 1$ we define
the Poisson and Hilbert-Poisson transforms by
\[P^r_{\zz}a=\sum_{n\in\zz} r^{|n|}z^n a_n,\qquad
H^r_{\zz}a=\sum_{n\in\zz} -(\sgn n)gr^{|n|}z^n a_n.\]
If $a=\sum_{n\in\zzh} z^n a_n $ then for  $0\leq r\leq 1$ we define
\[P^r_{\zzh}a=\sum_{n\in\zzh} r^{|n|}z^n a_n,\qquad
H^r_{\zzh}a=\sum_{n\in\zzh} -(\sgn n)gr^{|n|}z^n a_n.\]

In the  case $r=1$ these operations specialize as follows:
 $P^1_{\zzi}$ yields the identity, and  $H^1_{\zzi}$ yields the Hilbert transform.
\end{point}

\section{Some decompositions of matrices and related material}
\begin{point}
In case of $\mathsf L_{\zz}^{g,x_1}[u]$ the matrix was composed of two different parts.
We can generalize this construction as follows.
We define the matrices
\[\overleftarrow{\mathsf H}^{g}_{\zz,\zzh}
=\sum_{n\in\zzh^-} g\mathbf e_{n-\frac12,n}-\sum_{n\in\zzh^+} g\mathbf e_{n+\frac12,n},\qquad
\overrightarrow{\mathsf H}^{g}_{\zzh,\zz}=(\overleftarrow{\mathsf H}^{g}_{\zz,\zzh})^\top,\]
\[\overleftarrow{\mathsf H}^g_{\zzh,\zz}
=\sum_{n\in\zz^-} g\mathbf e_{n-\frac12,n}+\tfrac{\sqrt2}2g\mathbf e_{-\frac12,0}
-\tfrac{\sqrt2}2g\mathbf e_{\frac12,0}-\sum_{n\in\zz^+} g\mathbf e_{n+\frac12,n},\qquad
\overrightarrow{\mathsf H}^g_{\zz,\zzh}=(\overleftarrow{\mathsf H}^g_{\zzh,\zz})^\top;\]
and we let
$\mathbf f^g_{\zz}=\mathbf e_{0}, \mathbf f^g_{\zzh}=
\tfrac{\sqrt2}2\mathbf e_{-\frac12}+\tfrac{\sqrt2}2\mathbf e_{\frac12},
\bar{\mathbf f}^g_{\zz}=\mathbf e_{0}^\top, \bar{\mathbf f}^g_{\zzh}=
\tfrac{\sqrt2}2\mathbf e_{-\frac12}^\top+\tfrac{\sqrt2}2\mathbf e_{\frac12}^\top.$

If $\mathbf A$ is a $\zzi'\times \zzj'$ matrix, $\mathbf b$ is a $\zzi'$-column matrix,
$\mathbf c$ is a $\zzj'$-row matrix, $d$ is a scalar matrix, then we define
\[\left\{\begin{matrix}\mathbf A & \mathbf b\\\mathbf c& d \end{matrix}\right\}^g_{\zzi,\zzj}=
\overleftarrow{\mathsf H}^g_{\zzi,\zzi'}\mathbf A\overrightarrow{\mathsf H}^{-g}_{\zzj',\zzj}+
\overleftarrow{\mathsf H}^g_{\zzi,\zzi'}\mathbf b\otimes \bar{\mathbf f}^g_{\zzj}+
\mathbf f^g_{\zzi}\otimes\mathbf c\overrightarrow{\mathsf H}^{-g}_{\zzj',\zzj}+
d\mathbf f^g_{\zzi}\otimes \bar{\mathbf f}^g_{\zzj}.\]
These symbols multiply as $2\times 2$ matrices.
Every $\zzi\times\zzj$ matrix $\mathbf C$ occurs so:
\[\mathbf C=
\left\{\begin{matrix}
\overrightarrow{\mathsf H}^{-g}_{\zzi',\zzi}\mathbf C \overleftarrow{\mathsf H}^{g}_{\zzj',\zzj} &
\overrightarrow{\mathsf H}^{-g}_{\zzi',\zzi}\mathbf C \mathbf f^g_{\zzj} \\
\bar{\mathbf f}^g_{\zzi}\mathbf C \overleftarrow{\mathsf H}^{g}_{\zzj',\zzj} &
\bar{\mathbf f}^g_{\zzi}\mathbf C\mathbf f^g_{\zzj}
\end{matrix}\right\}^g_{\zzi,\zzj}.\]

There are variants of this construction but without $\frac{\sqrt2}2$.
Let $\mathsf X^g_{\zz}=\mathbf 1_{\zz},$ and
\[\mathsf X^g_{\zzh}=\mathbf 1_{\zz\setminus\{-\frac12,\frac12\}}+
\frac{\sqrt2}2\frac{(\mathbf e_{-\frac12}-\mathbf e_{\frac12})
\otimes(\mathbf e_{-\frac12}^\top-\mathbf e_{\frac12}^\top)}2+
\sqrt2\frac{(\mathbf e_{-\frac12}+\mathbf e_{\frac12})
\otimes(\mathbf e_{-\frac12}^\top+\mathbf e_{\frac12}^\top)}2;\]
and we define
\[\overleftarrow{\mathsf H}^{R,g}_{\zzi',\zzi}=\mathsf X^g_{\zzi'}\overleftarrow{\mathsf H}^{g}_{\zzi',\zzi},\quad
\mathbf f_{\zzi}^{R,g}=\mathsf X^g_{\zzi}\mathbf f^g_{\zzi},\quad
\overrightarrow{\mathsf H}^{R,g}_{\zzi,\zzi'}
=\overrightarrow{\mathsf H}^{g}_{\zzi,\zzi'}(\mathsf X^g_{\zzi'})^{-1},\quad
\bar{\mathbf f}_{\zzi}^{R,g}=\bar{\mathbf f}^g_{\zzi}(\mathsf X^g_{\zzi})^{-1},\]
\[\overleftarrow{\mathsf H}^{L,g}_{\zzi',\zzi}=(\mathsf X^g_{\zzi'})^{-1}
\overleftarrow{\mathsf H}^{g}_{\zzi',\zzi},\quad
\mathbf f_{\zzi}^{L,g}=(\mathsf X^g_{\zzi})^{-1}\mathbf f^g_{\zzi},\quad
\overrightarrow{\mathsf H}^{L,g}_{\zzi,\zzi'}
=\overrightarrow{\mathsf H}^{g}_{\zzi,\zzi'}\mathsf X^g_{\zzi'},\quad
\bar{\mathbf f}_{\zzi}^{L,g}=\bar{\mathbf f}^g_{\zzi}\mathsf X^g_{\zzi}.\]
Practically this  means that every term $\frac{\sqrt2}2$ is replaced by $1$ or $\frac12$.

Then we can define
$\left\{\begin{matrix}\mathbf A & \mathbf b\\\mathbf c& d \end{matrix}\right\}^{L,g}_{\zzi,\zzj}$
and  $\left\{\begin{matrix}\mathbf A & \mathbf b\\\mathbf c& d \end{matrix}\right\}^{R,g}_{\zzi,\zzj}$
analogously.
In fact, the relationship between the basic version and these is very simple:
\[\left\{\begin{matrix}\mathbf A & \mathbf b\\\mathbf c& d \end{matrix}\right\}^{R,g}_{\zzi,\zzj}\!\!=
\mathsf X^g_{\zzi}
\left\{\begin{matrix}\mathbf A & \mathbf b\\\mathbf c& d \end{matrix}\right\}^{g}_{\zzi,\zzj}\!\!
(\mathsf X^g_{\zzi})^{-1},\!\!\qquad\!\!
\left\{\begin{matrix}\mathbf A & \mathbf b\\\mathbf c& d \end{matrix}\right\}^{L,g}_{\zzi,\zzj}\!\!=
(\mathsf X^g_{\zzi})^{-1}
\left\{\begin{matrix}\mathbf A & \mathbf b\\\mathbf c& d \end{matrix}\right\}^{g}_{\zzi,\zzj}\!\!
\mathsf X^g_{\zzi}.
\]

The relationship to our original  motivating example is that
\[\mathsf L^{g,x_1}_{\zz}[u]
=\left\{\begin{matrix}\mathsf L^{g,x_1}_{\zzh} &\\& u \end{matrix}\right\}^{g}_{\zz}
=\left\{\begin{matrix}\mathsf L^{g,x_1}_{\zzh} &\\& u \end{matrix}\right\}^{R,g}_{\zz}
=\left\{\begin{matrix}\mathsf L^{g,x_1}_{\zzh} &\\& u \end{matrix}\right\}^{L,g}_{\zz}.\]
In general,
$\left\{\begin{matrix}\mathbf A&\mathbf b\\\mathbf c& d \end{matrix}\right\}^{X,g}_{\zz}$
and $\left\{\begin{matrix} d_1\mathbf 1_{\zzi'}&\\& d_2 \end{matrix}\right\}^{X,g}_{\zzi}$
do not depend on the choice of $X$, which can be $\emptis$ (the empty symbol), $R$ or $L$.
\end{point}

\begin{point} \textit{Difference functions.}
For $a\in\mathfrak B^0$ we consider the difference functions
\[\Delta_{\zz,\zz}^{x_1}a(x)=\frac12 \frac{(a(x)-a(x_1))\cos\frac{x-x_1}{2}}{\sin\frac{x-x_1}{2}},\]
\[\Delta_{\zzh,\zz}^{x_1}a(x)=\frac12 \frac{a(x)-a(x_1)}{\sin\frac{x-x_1}{2}}, \]
and for  $a\in\mathfrak B^1$
\[\Delta_{\zz,\zzh}^{x_1}a(x)=\frac12 \frac{a(x)-a(x_1)\cos\frac{x-x_1}{2}}{\sin\frac{x-x_1}{2}},\]
\[\Delta_{\zzh,\zzh}^{x_1}a(x)=\frac12 \frac{a(x)\cos\frac{x-x_1}{2}-a(x_1)}{\sin\frac{x-x_1}{2}}.\]
\end{point}
\begin{lemma}
For $n\in\zz^+$ it yields
\begin{multline}
\Delta^{x_1}_{\zz,\zz} (z^{-n})=
-g\left(\tfrac12z^{-n}+z_1^{-1}z^{-n+1}+\ldots+z_1^{-n+1}z^{-1}+\tfrac12z_1^{-n}\right),
\\ \Delta^{x_1}_{\zz,\zz} (1)= 0,  \\
\Delta^{x_1}_{\zz,\zz} (z^{n})=g\left(\tfrac12z_1^n+z_1^{n-1}z+\ldots+z_1z^{n-1}+\tfrac12z^n\right),
\notag\end{multline}
\begin{multline}
\Delta^{x_1}_{\zzh,\zz} (z^{-n})=
-g\left(z_1^{-\frac12}z^{-n+\frac12}+\ldots+z_1^{-n+\frac12}z^{-\frac12}\right),
\\\Delta^{x_1}_{\zzh,\zz} (1)= 0,  \\
\Delta^{x_1}_{\zzh,\zz} (z^{n})=g\left(z_1^{n-\frac12}z^{\frac12}+\ldots+z_1^{\frac12}z^{n-\frac12}\right);
\notag\end{multline}
and for $n\in\zzh^+$ it yields
\begin{multline}
\Delta^{x_1}_{\zz,\zzh} (z^{-n})=
-g\left(z_1^{-\frac12}z^{-n+\frac12}+\ldots +z_1^{-n+1}z^{-1}+\tfrac12z_1^{-n}\right),
\\
\Delta^{x_1}_{\zz,\zzh} (z^n)=
g\left(\tfrac12z_1^{n}+ z_1^{n-1}z +\ldots +z_1^{\frac12}z^{n-\frac12}\right),
\notag\end{multline}
\begin{multline}
\Delta^{x_1}_{\zzh,\zzh} (z^{-n})=
-g\left({\tfrac12}z^{-n}+ z_1^{-1} z^{-n+1} +\ldots+z_1^{-n+\frac12}z^{-\frac12}\right),
\\
\Delta^{x_1}_{\zzh,\zzh} (z^n)=
g\left(z_1^{n-\frac12}z^{\frac12}+\ldots+ z_1 z^{n-1} +{\tfrac12}z^{n}\right).
\notag\end{multline}
\qed\end{lemma}

\begin{lemma}It yields
\[\int\Delta^{x_1}_{\zz,\zzl}a(x)\,\frac{dx}{2\pi}=-\frac12 H_{\zzl}a(x_1), \qquad
\int\Delta^{x_1}_{\zzl,\zzl}a(x)\,\frac{dx_1}{2\pi}=-\frac12 H_{\zzl}a(x).\]
\qed\end{lemma}
\begin{point}
We define the critical part of the difference functions as
\[\widehat\Delta^{x_1}_{\zzi,\zzl}a (x)=
-\tfrac12\delta_{i,0} H_{\zzl}a(x_1) -\tfrac12\delta_{j,0} H_{\zzl}a(x).  \]
\end{point}

\begin{point}
\textit{Multiplicative difference functions.}
If $a$ is multiplicatively invertible then it is also natural to consider
multiplicative difference functions:
For $a\in\mathfrak B^0$
\[\Xi_{\zz,\zz}a(x,x_1)=\frac12 \frac{(a(x)a(x_1)^{-1}-1)\cos\frac{x-x_1}{2}}{\sin\frac{x-x_1}{2}},\]
\[\Xi_{\zzh,\zz}a(x,x_1)=\frac12 \frac{a(x)a(x_1)^{-1}-1}{\sin\frac{x-x_1}{2}}, \]
and for  $a\in\mathfrak B^1$
\[\Xi_{\zz,\zzh}a(x,x_1)=\frac12 \frac{a(x)a(x_1)^{-1}-\cos\frac{x-x_1}{2}}{\sin\frac{x-x_1}{2}},\]
\[\Xi_{\zzh,\zzh}a(x,x_1)=\frac12 \frac{a(x)a(x_1)^{-1}\cos\frac{x-x_1}{2}-1}{\sin\frac{x-x_1}{2}}.\]
\begin{lemma}
\[\int\Xi_{\zz,\zzl}a(x,x_1)\,\frac{dx}{2\pi}=-\frac12 H_{\zzl}a(x_1)\cdot a(x_1)^{-1},  \]
\[\int\Xi_{\zz,\zzl}a(x_1,x)\,\frac{dx}{2\pi}=\frac12 a(x_1)\cdot H_{\zzl}(a^{-1})(x_1).  \]
\qed\end{lemma}
\begin{point}
We can define normalized multiplicative difference functions by
\[\Xi^L_{\zzi,\zzl}a (x,x_1)=
\Xi_{\zzi,\zzl}a (x,x_1)+\tfrac12\delta_{i,0}
 H_{\zzl}a(x_1)a(x_1)^{-1} +\tfrac12\delta_{j,0}H_{\zzl}a(x)a(x_1)^{-1},  \]
\[\Xi^R_{\zzi,\zzl}a(x_1,x)=\Xi_{\zzi,\zzl}a(x_1,x)
)-\tfrac12\delta_{i,0}a(x_1)H_{\zzl}a^{-1}(x_1) -\tfrac12\delta_{j,0}a(x_1)H_{\zzl}a^{-1}(x).  \]
Then
\[\int\Xi_{\zz,\zzl}^La(x,x_1)\,\frac{dx}{2\pi}=0,\qquad\text{and}\qquad
\int\Xi_{\zz,\zzl}^Ra(x_1,x)\,\frac{dx}{2\pi}=0.\]
\end{point}
\begin{point}
The following table might help to keep track the various periodicity properties:
\[\begin{array}{|c|c|c|c|c|c|}
\hline l & i & j &  a(x)  \text{ in }  x &
\Delta_{\zzi,\zzl}^{x_1}a(x) \text{ in } x  & \Delta_{\zzi,\zzl}^{x_1}a(x) \text{ in } x_1 \\
 &  &  &  &
\Delta_{\zzj,\zzl}^{x_1}a(x) \text{ in } x_1  & \Delta_{\zzj,\zzl}^{x_1}a(x) \text{ in } x \\
& & & & \Xi_{\zzi,\zzl}(x,x_1) \text{ in }x\text{ or }x_1 &
\Xi_{\zzj,\zzl}(x_1,x) \text{ in }x\text{ or }x_1\\
\hline
0&0&0&per&per&per\\
0&1&1&per&skew&skew\\
1&0&1&skew&per&skew\\
1&1&0&skew&skew&per\\\hline
\end{array}\]
\end{point}

\end{point}

\begin{lemma}\label{lem:precorr}
It yields
\[\left\{\begin{matrix}\mathsf U_{\zzi',\zzj'}^g(a) & \\ \bar{\mathsf S}^g_{\zzj'}(\Delta^{x_1}_{\zzj',\zzl}a)
& a(x_1) \end{matrix}\right\}^{X,g}_{\zzi,\zzj}=\]
\[=\left\{\begin{matrix} \mathbf 1_{\zzi'}& \\ \bar{\mathsf S}^g_{\zzi'}(\Xi_{\zzi',\zzl}a(x_1,x))
& 1\end{matrix}\right\}^{X,g}_{\zzi}
\left\{\begin{matrix}\mathsf U_{\zzi',\zzj'}^g(a) & \\& a(x_1) \end{matrix}\right\}^{X,g}_{\zzi,\zzj}=\]
\[=\left\{\begin{matrix}\mathsf U_{\zzi',\zzj'}^g(a) & \\& a(x_1) \end{matrix}\right\}^{X,g}_{\zzi,\zzj}
\left\{\begin{matrix} \mathbf 1_{\zzj'}& \\ -\bar{\mathsf S}^g_{\zzj'}(\Xi_{\zzj',\zzl}(a^{-1})(x_1,x))
& 1\end{matrix}\right\}^{X,g}_{\zzj}=\]
\[=\left\{\begin{matrix} \mathbf 1_{\zzi'}& \\ \bar{\mathsf S}^g_{\zzi'}(\Xi^L_{\zzi',\zzl}a(x_1,x))
& 1\end{matrix}\right\}^{X,g}_{\zzi}
\left\{\begin{matrix}\mathsf U_{\zzi',\zzj'}^g(a) & \\
\bar{\mathsf S}_{\zzj'}^g(\widehat\Delta^{x_1}_{\zzj',\zzl}a(x))
& a(x_1) \end{matrix}\right\}^{X,g}_{\zzi,\zzj}=\]
\[=\left\{\begin{matrix}\mathsf U_{\zzi',\zzj'}^g(a) & \\
\bar{\mathsf S}^g_{\zzj'}(\widehat\Delta^{x_1}_{\zzj',\zzl}a(x))
& a(x_1) \end{matrix}\right\}^{X,g}_{\zzi,\zzj}
\left\{\begin{matrix} \mathbf 1_{\zzj'}& \\ -\bar{\mathsf S}^g_{\zzj'}(\Xi^R_{\zzj',\zzl}(a^{-1})(x_1,x))
& 1\end{matrix}\right\}^{X,g}_{\zzj}=\]
\[=\left\{\begin{matrix}
\mathbf 1_{\zzi'}&\\\frac12\bar{\mathsf h}_{\zzi'}^{g,x_1}&1 \end{matrix}\right\}^{X,g}_{\zzi}
\left\{\begin{matrix}
\mathsf U_{\zzi',\zzj'}^g(a) & \\& a(x_1) \end{matrix}\right\}^{X,g}_{\zzi,\zzj}
\left\{\begin{matrix}
\mathbf 1_{\zzj'}&\\-\frac12\bar{\mathsf h}_{\zzj'}^{g,x_1}&1 \end{matrix}\right\}^{X,g}_{\zzj}.\]
and
\[\left\{\begin{matrix}\mathsf U_{\zzi',\zzj'}^g(a) &{\mathsf S}^g_{\zzi'}(\Delta^{x_1}_{\zzi',\zzl}a)
\\& a(x_1) \end{matrix}\right\}^{X,g}_{\zzi,\zzj}=\]
\[=\left\{\begin{matrix}\mathbf 1_{\zzi'} &{\mathsf S}^g_{\zzi'}(\Xi_{\zzi',\zzl}a(x,x_1))
\\& 1 \end{matrix}\right\}^{X,g}_{\zzi}
\left\{\begin{matrix}\mathsf U_{\zzi',\zzj'}^g(a) &
\\& a(x_1) \end{matrix}\right\}^{X,g}_{\zzi,\zzj}=\]
\[=\left\{\begin{matrix}\mathsf U_{\zzi',\zzj'}^g(a) &\\& a(x_1) \end{matrix}\right\}^{X,g}_{\zzi,\zzj}
\left\{\begin{matrix}\mathbf 1_{\zzj'} &-{\mathsf S}^g_{\zzj'}(\Xi_{\zzj',\zzl}(a^{-1})(x,x_1))
\\& 1 \end{matrix}\right\}^{X,g}_{\zzj}=\]
\[=\left\{\begin{matrix}\mathbf 1_{\zzi'} &{\mathsf S}^g_{\zzi'}(\Xi^L_{\zzi'}a(x,x_1))
\\& 1 \end{matrix}\right\}^{X,g}_{\zzi}
\left\{\begin{matrix}\mathsf U_{\zzi',\zzj'}^g(a) &
\mathsf S^g_{\zzi'}(\widehat\Delta^{x_1}_{\zzi',\zzl}a(x))
\\& a(x_1) \end{matrix}\right\}^{X,g}_{\zzi,\zzj}=\]
\[=\left\{\begin{matrix}\mathsf U_{\zzi',\zzj'}^g(a) &
\mathsf S^g_{\zzi'}(\widehat\Delta^{x_1}_{\zzi',\zzl}a(x))
\\& a(x_1) \end{matrix}\right\}^{X,g}_{\zzi,\zzj}
\left\{\begin{matrix}\mathbf 1_{\zzj'} &-{\mathsf S}^g_{\zzj'}(\Xi^R_{\zzj',\zzl}(a^{-1})(x,x_1))
\\& 1 \end{matrix}\right\}^{X,g}_{\zzj}=\]
\[=\left\{\begin{matrix}\mathbf 1_{\zzi'}&-\frac12{\mathsf h}_{\zzi'}^{g,x_1}\\&1 \end{matrix}\right\}^{X,g}_{\zzi}
\left\{\begin{matrix}\mathsf U_{\zzi',\zzj'}^g(a) &\\& a(x_1) \end{matrix}\right\}^{X,g}_{\zzi,\zzj}
\left\{\begin{matrix}\mathbf 1_{\zzj'}&\frac12{\mathsf h}_{\zzj'}^{g,x_1}\\&1 \end{matrix}\right\}^{X,g}_{\zzj}.\]
\qed\end{lemma}
In fact, the matrices involving row or column matrices of principal value distributions of simple poles
are not rapidly decreasing but one can compute with them legally in an appropriate topology.

\begin{lemma}\label{lem:excorr}
The entries of the matrices
\[\mathsf L^{g,\frac{x_1}2}_{\zzi}
\left\{\begin{matrix}\mathsf U_{\zzi',\zzj'}^g(a) & \mathsf S^g_{\zzi'}(
\widehat\Delta_{\zzi',\zzl}^{x_1}a(x)  )\\& a(x_1) \end{matrix}\right\}^{L,g}_{\zzi,\zzj}
\mathsf L^{g,\frac{x_1}2}_{\zzj}\]
and
\[\mathsf L^{g,\frac{x_1}2}_{\zzi}
\left\{\begin{matrix}\mathsf U_{\zzi',\zzj'}^g(a) &\\ \bar{\mathsf S}^g_{\zzj'}(
\widehat\Delta_{\zzj',\zzl}^{x_1}a(x)  )& a(x_1) \end{matrix}\right\}^{R,g}_{\zzi,\zzj}
\mathsf L^{g,\frac{x_1}2}_{\zzj}\]
are possibly infinite sums of terms of shape $\pm[a_{p}]^\pm_gz_1^{q}$.
\qed\end{lemma}
The lemma above is valid without the ``critical'' mark $\widehat{~~}$ but in that case the
statement does not require particular checking, as we will see.

\section{Some matrices built up from elementary rotations}
In this section we define some natural matrices.
\begin{point}
For $n\in\frac12\zz$ and a ``scalar'' $A$
we let $\mathsf I_n(A)$ be the matrix which is $\mathbf 1_{(\zz+n)}$ but the
coefficient of $\mathbf e_{n,n}$ is replaced by $A$, ie.
\[\mathsf I_n(A)=\mathbf 1_{(\zz+n)\setminus\{n\}}+A\mathbf e_{n,n}.\]

For $n\in\frac12\zz$ and a $2\times2$ matrix $A$
we let $\notI_n(A)$ be the matrix which is $\mathbf 1_{(\zzh+n)}$ but the
$\{n-\frac12,n+\frac12\}\times\{n-\frac12,n+\frac12\}$ block is replaced by $A$.
More precisely, we set
\[\notI_n\!\left(\begin{bmatrix}a&b\\c&d\end{bmatrix}\right)
=\mathbf 1_{(\zzh+n)\setminus\{n-\frac12,n+\frac12\}}+
a\,\mathbf e_{n-\frac12,n-\frac12}+b\,\mathbf e_{n-\frac12,n+\frac12}
+c\,\mathbf e_{n+\frac12,n-\frac12}+d\,\mathbf e_{n+\frac12,n+\frac12}.\]
\end{point}
\begin{point}We define
\[\mathsf A_{n}(\theta)
=\notI_n\!\left(\begin{bmatrix} \cos\theta & \sin\theta\\ -\sin\theta& \cos\theta \end{bmatrix}\right)
,\qquad
\daa_{n}(\theta)=\notI_n\!\left(
\begin{bmatrix} \cos\frac\theta2 & \sin\frac\theta2\\ -\sin\frac\theta2& \cos\frac\theta2 \end{bmatrix}\right),\]
where $\theta\in\mathbb R$.
Similarly, we can define
\[\mathsf A^R_{n}(\theta)=
\notI_n\!\left(\begin{bmatrix} 1 & \sin\theta\\ -\sin\theta& 1-\sin^2\theta \end{bmatrix}\right),
\qquad\daa^R_{n}(\theta)=
\notI_n\!\left(\begin{bmatrix} 1 & \\ -\sin\theta& 1\end{bmatrix}\right),\]
\[\mathsf A^L_{n}(\theta)=
\notI_n\!\left(\begin{bmatrix} 1-\sin^2\theta & \sin\theta\\ -\sin\theta& 1 \end{bmatrix}\right),
\qquad\daa^L_{n}(\theta)=
\notI_n\!\left(\begin{bmatrix} 1 & \sin\theta\\ & 1 \end{bmatrix}\right).\]
In this case we can use natural parametrization by $t=\sin\theta$.
We present yet another version, which is less nice than the previous ones
but it may illustrate the structure of these matrices in question a bit further.
We define
\[\mathsf A^{\frac R3}_{n}(\theta)=
\notI_n\!\left(\begin{bmatrix}\cos\theta\cos\frac\theta3&\sin\theta\\
-\sin\theta&\frac{\cos\theta}{\cos\frac\theta3}\end{bmatrix}\right),\qquad
\daa^{\frac R3}_{n}(\theta)=
\notI_n\!\left(\begin{bmatrix}\cos\frac{2\theta}3&\sin\frac{\theta}3\\-2\sin\frac{\theta}3&1\end{bmatrix}\right),\]
\[\mathsf A^{\frac L3}_{n}(\theta)=
\notI_n\!\left(\begin{bmatrix}\frac{\cos\theta}{\cos\frac\theta3}&\sin\theta\\
-\sin\theta&\cos\theta\cos\frac\theta3\end{bmatrix}\right),\qquad
\daa^{\frac L3}_{n}(\theta)=
\notI_n\!\left(\begin{bmatrix}1&2\sin\frac{\theta}3\\-\sin\frac{\theta}3&\cos\frac{2\theta}3\end{bmatrix}\right).\]
In this case we can use natural parametrization by $\tau=2\sin\frac\theta3$.
Then
\[\sin\frac\theta3=\frac\tau2,\qquad
\cos\frac{2\theta}3=1-\frac{\tau^2}2,\qquad
\frac{\cos\theta}{\cos\frac\theta3}=2\cos\frac{2\theta}3-1=1-\tau^2,\]
\[\cos\theta\cos\frac\theta3=\frac12\left(\cos\frac{2\theta}3+\cos\frac{4\theta}3\right)=
\frac{(1-\tau^2)(4-\tau^2)}4,\qquad
\sin\theta=\frac{\tau(3-\tau^2)}2;\]
so the entries of the replacement blocks are simple expressions of $\tau$.
\end{point}
\begin{point}
We will use the symbol $X$ for any of $R,\frac R3,\emptis$ (the empty symbol), $\frac L3, L$.
Then  the symbol  $-X$ will be used for  $L,\frac L3,\emptis,\frac R3, R$, respectively.

The matrices $\mathsf A^X_{n}(\theta)$ are to be considered as elementary rotations, and
the matrices $\daa^X_{n}(\theta)$ are to be considered as elementary half-rotations.
It yields
\[\mathsf A^X_n (\theta)=\daa^X_n(\theta)\daa^{-X}_n(-\theta)^{-1},\quad
\daa^X_n (\theta)^\top=\daa^{-X}_n (\theta)^{-1},\quad
\fli_{\zz_n'} \daa^X_n(\theta)\fli_{\zz_n'}= \daa^{-X}_{-n}(-\theta);\]
hence
\[\mathsf A^X_n(\theta)=\mathsf A^{-X}_n(-\theta)^{-1},\qquad
\mathsf A^X_n(\theta)^\top=\mathsf A^X_n(-\theta),\qquad
\fli_{\zz_n'} \mathsf A^X_n(\theta)\fli_{\zz_n'}= \mathsf A^{X}_{-n}(\theta)^{-1};\]
although for particular choices of $X$ more special equalities hold.
\end{point}
In the rest of the section we assume that
$\theta\notin\pi\zzh$, ie.~that $\sin\theta\neq\pm1$.
\begin{point}
In what follows, it is useful to consider the function
\[\omega(x,\theta)\equiv \omega(z,z^{-1},\theta)
=(1-z\sin\theta)(1-z^{-1}\sin\theta)=1-2\cos x\sin\theta+\sin^2\theta;\]
even if only as an abbreviation. It is a positive real function.
In particular, we can take various powers of it without the problem of ambiguities.
\end{point}
\begin{point}
Let us consider the matrices
$\widetilde{\mathsf A}^{R,R,g}_{\zz}(\theta)=\mathbf 1_{\zz}$,
\[\widetilde{\mathsf A}^{\frac R3,R,g}_{\zz}(\theta)
=\mathsf I_0\!\left(\frac{\cos\theta}{\cos\frac\theta3}\right),\qquad
\widetilde{\mathsf A}^{\emptis,R,g}_{\zz}(\theta)
=\mathsf I_0\!\left(\cos\theta\right),\]
\[\widetilde{\mathsf A}^{\frac L3,R,g}_{\zz}(\theta)
=\mathsf I_0\!\left(\cos\theta\cos\frac\theta3\right),\qquad
\widetilde{\mathsf A}^{L,R,g}_{\zz}(\theta)
=\mathsf I_0\!\left(\cos^2 \theta\right);\]
and
$\widetilde{\mathsf A}^{R,R,g}_{\zzh}(\theta)=\mathbf 1_{\zzh}$,
\[\widetilde{\mathsf A}^{\frac R3,R,g}_{\zzh}(\theta)=
\notI_0\!\left(\begin{bmatrix} \cos\tfrac{2\theta}3 & -\sin\tfrac\theta3\\
-\sin\tfrac\theta3&\cos\tfrac{2\theta}3 \end{bmatrix}\right),\qquad
\widetilde{\mathsf A}^{\emptis,R,g}_{\zzh}(\theta)=
\notI_0\!\left(\begin{bmatrix} \cos\tfrac{\theta}2 & -\sin\tfrac\theta2\\
-\sin\tfrac\theta2&\cos\tfrac{\theta}2 \end{bmatrix}\right),\]
\[\widetilde{\mathsf A}^{\frac L3,R,g}_{\zzh}(\theta)=
\notI_0\!\left(\begin{bmatrix} 1 & -2\sin\tfrac\theta3\\-2\sin\tfrac\theta3&1 \end{bmatrix}\right),\qquad
\widetilde{\mathsf A}^{L,R,g}_{\zzh}(\theta)=
\notI_0\!\left(\begin{bmatrix} 1 & -\sin\theta\\-\sin\theta&1 \end{bmatrix}\right).\]
One can extend this notation so that
$\widetilde{\mathsf A}^{X,Y,g}_{\zzi}(\theta)=\widetilde{\mathsf A}^{X,Z,g}_{\zzi}(\theta)
\widetilde{\mathsf A}^{Z,Y,g}_{\zzi}(\theta)$
would hold.
\end{point}
\begin{point}
Consider the matrices
\begin{multline}
\mathsf F^{R,g}_{\zz}(\theta)=
\mathbf 1_{\zz^-}\mathsf U^g_{\zz}\left(1-z^{-1}\sin\theta\right)+\mathbf e_{0,0}
+\mathbf 1_{\zz^+}\mathsf U^g_{\zz}\left(1-z\sin\theta \right)=\\
=\mathbf 1_{\zz}-\mathsf U^g_{\zz}(\cos x)\sin\theta+\mathsf H^g_{\zz}
\mathsf U^g_{\zz}(\sin x)\sin\theta+\mathbf e_{0,0}\mathsf U^g_{\zz}(\cos x)\sin\theta=\\=
\mathbf 1_{\zz}-\mathsf U^g_{\zz}(\cos x)\sin\theta+
\mathsf U^g_{\zz}(\sin x)\mathsf H^g_{\zz}\sin\theta
-\mathbf e_{0,0}\mathsf U^g_{\zz}(\cos x)\sin\theta,
\notag\end{multline}
\begin{multline}
\mathsf F^{R,g}_{\zzh}(\theta)=
\mathbf 1_{\zzh^{--}}\mathsf U^g_{\zz}\left(1-z^{-1}\sin\theta\right)
+\mathbf 1_{\{-\frac12,\frac12\}}
+\mathbf 1_{\zzh^{++}}\mathsf U^g_{\zz}\left(1-z\sin\theta \right) =\\
=\mathsf U^g_{\zz}\left(1-z^{-1}\sin\theta\right)\mathbf 1_{\zzh^{-}}
+\mathsf U^g_{\zz}\left(1-z\sin\theta \right)\mathbf 1_{\zzh^{+}}=\\
=\mathbf 1_{\zz}-\mathsf U_{\zz}^g(\cos x)\sin\theta+
\mathsf U^g_{\zzh}(\sin x)\mathsf H^g_{\zzh}\sin\theta;
\notag\end{multline}
and, furthermore, we define the matrices
\[\mathsf F^{X,g}_{\zzi}(\theta)
=\widetilde{\mathsf A}^{X,R,g}_{\zzi}(\theta)\mathsf F^{R,g}_{\zzi}(\theta)
=\widetilde{\mathsf A}^{X,R,g}_{\zzi}(\theta)+\mathsf F^{R,g}_{\zzi}(\theta)-\mathbf 1_{\zzi}.\]
In particular, it yields
\begin{multline}
\mathsf F^g_{\zz}(\theta)
=\mathbf 1_{\zz^-}\mathsf U^g_{\zz}\left(1-z^{-1}\sin\theta\right)+\mathbf e_{0,0}\cos\theta
+\mathbf 1_{\zz^+}\mathsf U^g_{\zz}\left(1-z\sin\theta \right),\hfill
\notag\end{multline}
\begin{multline}
\mathsf F^g_{\zzh}(\theta)=
\mathbf 1_{\zzh^{--}}\mathsf U^g_{\zzh}
\left(1-z^{-1}\sin\theta\right)+\mathbf e_{-\frac12,-\frac12}
\mathsf U^g_{\zzh}\left(\cos\tfrac\theta2-z^{-1}\sin\tfrac\theta2\right)+\\
+\mathbf e_{\frac12,\frac12}\mathsf U^g_{\zzh}\left(\cos\tfrac\theta2-z\sin\tfrac\theta2\right)
+\mathbf 1_{\zzh^{++}}\mathsf U^g_{\zzh}
\left(1-z\sin\theta\right),
\notag\end{multline}
\begin{multline}
\mathsf F^{L,g}_{\zz}(\theta)
=\mathbf 1_{\zz^-}\mathsf U^g_{\zz}\left(1-z^{-1}\sin\theta\right)
+\mathbf e_{0,0}(1-\sin^2\theta)
+\mathbf 1_{\zz^+}\mathsf U^g_{\zz}\left(1-z\sin\theta \right),\hfill
\notag\end{multline}
\begin{multline}
\mathsf F^{L,g}_{\zzh}(\theta)
=\mathbf 1_{\zzh^{-}}\mathsf U^g_{\zzh}\left(1-z^{-1}\sin\theta\right)
+\mathbf 1_{\zzh^{+}}\mathsf U^g_{\zzh}\left(1-z\sin\theta\right)=\\
=\mathbf 1_{\zz}-\mathsf U_{\zz}^g(\cos x)\sin\theta+
\mathsf H^g_{\zzh}\mathsf U^g_{\zzh}(\sin x)\sin\theta.
\notag\end{multline}
We also define
\[\mathsf F'^{X,g}_{\zzi}(\theta)=\mathsf F^{-X,g}_{\zzi}(\theta)^\top.\]
\end{point}
\begin{lemma}It yields
\[\mathsf F^{X,g}_{\zzi}(\theta)\mathsf U_{\zzi}^g(\omega(x,\theta)^{-1})
\mathsf F'^{X,g}_{\zzi}(\theta)=\mathbf 1_{\zzi}.\]
\qed\end{lemma}
\begin{point}
For $c=-1,-\frac13,0,\frac13,1$ let $cL$ denote $R,\frac R3,\emptis,\frac L3,L$ respectively.
Let
\[\mathsf D^{cL,g}_{\zzi}(\theta)=
\mathsf F^{cL,g}_{\zzi}(\theta)\mathsf U_{\zzi}^g(\omega(x,\theta)^{-\frac{1+c}2}).\]

In particular, it yields
\begin{multline}
\mathsf D^{R,g}_{\zz}(\theta)=\mathbf 1_{\zz^-}\mathsf U^g_{\zz}
\left(1-z^{-1}\sin\theta\right)+\mathbf e_{0,0}+
\mathbf 1_{\zz^+}\mathsf U^g_{\zz}\left(1-z\sin\theta\right),\hfill
\notag\end{multline}
\begin{multline}
\mathsf D^{R,g}_{\zzh}(\theta)=\mathbf 1_{\zzh^{--}}
\mathsf U^g_{\zzh}\left(1-z^{-1}\sin\theta\right)+\mathbf 1_{\{-\frac12,\frac12\}}
+\mathbf 1_{\zzh^{++}}\mathsf U^g_{\zzh}\left(1-z\sin\theta\right),\hfill
\notag\end{multline}
\begin{multline}
\mathsf D^{g}_{\zz}(\theta)
=\mathbf 1_{\zz^-}\mathsf U^g_{\zz}\left(\sqrt{\frac{1-z^{-1}\sin\theta}{1-z\sin\theta}}\right)
+\mathbf e_{0,0}\mathsf U^g_{\zz}\left(\frac{\cos\theta}{\sqrt{\omega(z,z^{-1},\theta)}}\right)+\\
+\mathbf 1_{\zz^+}\mathsf U^g_{\zz}\left(\sqrt{\frac{1-z\sin\theta}{1-z^{-1}\sin\theta}}\right),
\notag\end{multline}
\begin{multline}
\mathsf D^{g}_{\zzh}(\theta)
=\mathbf 1_{\zzh^{--}}\mathsf U^g_{\zzh}
\left(\sqrt{\frac{1-z^{-1}\sin\theta}{1-z\sin\theta}}\right)+\mathbf e_{-\frac12,-\frac12}
\mathsf U^g_{\zzh}\left(\frac{\cos\frac\theta2-z^{-1}\sin\frac\theta2}{
\sqrt{\omega(z,z^{-1},\theta)}}\right)+\\
+\mathbf e_{\frac12,\frac12}\mathsf U^g_{\zzh}\left(\frac{\cos\frac\theta2-z\sin\frac\theta2}{
\sqrt{\omega(z,z^{-1},\theta)}}\right)
+\mathbf 1_{\zzh^{++}}\mathsf U^g_{\zzh}
\left(\sqrt{\frac{1-z\sin\theta}{1-z^{-1}\sin\theta}}\right),
\notag\end{multline}
\begin{multline}
\mathsf D^{L,g}_{\zz}(\theta)=\mathbf 1_{\zz^-}
\mathsf U^g_{\zz}\left(\frac{1}{1-z\sin\theta}\right)
+\mathbf e_{0,0}\mathsf U^g_{\zz}
\left(\frac{1-\sin^2\theta}{\omega(z,z^{-1},\theta)}\right)+
\mathbf 1_{\zz^+}\mathsf U^g_{\zz}\left(\frac{1}{1-z^{-1}\sin\theta}\right),\hfill
\notag\end{multline}
\begin{multline}
\mathsf D^{L,g}_{\zzh}(\theta)
=\mathbf 1_{\zzh^{-}}\mathsf U^g_{\zzh}\left(\frac{1}{1-z\sin\theta}\right)
+\mathbf 1_{\zzh^{+}}\mathsf U^g_{\zzh}\left(\frac{1}{1-z^{-1}\sin\theta}\right).\hfill
\notag\end{multline}
One can check that
\[\mathsf D^{cL,g}_{\zzi}(\theta)^{-1}=
\mathsf U_{\zzi}^g(\omega(x,\theta)^{\frac{1-c}2})\mathsf F'^{cL,g}_{\zzi}(\theta);\]
or, in different terms,
\[\mathsf D^{X,g}_{\zzi}(\theta)^{-1}=\mathsf D^{-X,g}_{\zzi}(\theta)^\top.\]

Furthermore, $\fli_{\zzi}$ commutes with
$\mathsf D^{X,g}_{\zzi}(\theta)$.
\end{point}
\begin{point}
The matrices $\mathsf D_{\zzi}^{X,g}(\theta)$ and $\mathsf A_n^X(\theta)$
are closely related in the skew-complex representation.
If we take sufficiently weak limits, entrywise, for example, then we find
\[\mathsf U_{\zz}^g\left(\frac{1-z\sin\theta}{1-z^{-1}\sin\theta}\right)
=\wlim_{n,m\rightarrow\infty}
\mathsf A^X_{-m-\frac12}(\theta)\ldots\mathsf A^X_{-\frac12}(\theta)\mathsf A^X_{\frac12}(\theta)
\ldots \mathsf A^X_{\frac12+n}(\theta)\]
and
\[\mathsf U_{\zzh}^g\left(\frac{1-z\sin\theta}{1-z^{-1}\sin\theta}\right)
=\wlim_{n,m\rightarrow\infty}
\mathsf A^X_{-m}(\theta)\ldots
\mathsf A^X_{0}(\theta)
\ldots \mathsf A^X_{n}(\theta);\]
where $n,m\in\mathbb Z$.
Similarly, it yields
\[\mathsf F^{X,g}_{\zz}(\theta)\mathsf U_{\zz}^g\left(\frac{1}{1-z^{-1}\sin\theta}\right)=
\wlim_{n\rightarrow\infty}\mathsf A^X_{\frac12}(\theta)\ldots\mathsf A^X_{n+\frac12}(\theta)\]
and
\[\mathsf F^{X,g}_{\zzh}(\theta)\mathsf U_{\zzh}^g\left(\frac{1}{1-z^{-1}\sin\theta}\right)
=\wlim_{n\rightarrow\infty}\daa^{-X}_{0}(-\theta)^{-1}\mathsf A^X_{1}(\theta)\ldots \mathsf A^X_{n}(\theta).\]

Hence, by taking $\sqrt{A\cdot B}= B\sqrt{B^{-1}A}$, we obtain
\begin{multline}
\mathsf F^{X,g}_{\zz}(\theta)\mathsf U^{g}_{\zz}(\omega(x,\theta)^{-\frac12})=\\
=\sqrt{
\left(\wlim_{m\rightarrow\infty}\mathsf A^X_{-\frac12}(\theta)^{-1}\ldots
\mathsf A^X_{-m-\frac12}(\theta)^{-1}\right)
\cdot\left(\wlim_{n\rightarrow\infty}\mathsf A^X_{\frac12}(\theta)\ldots
\mathsf A^X_{n+\frac12}(\theta)\right)},
\notag\end{multline}
and
\begin{multline}
\mathsf F^{X,g}_{\zzh}(\theta)\mathsf U^{g}_{\zzh}(\omega(x,\theta)^{-\frac12})
=\\=\sqrt{\left(\wlim_{m\rightarrow\infty}\daa^{X}_{0}(\theta)^{-1}
\mathsf A^X_{-1}(\theta)^{-1}\ldots\mathsf A^X_{-m}(\theta)^{-1}\right)\cdot
\left(\wlim_{n\rightarrow\infty}\daa^{-X}_{0}(\theta)^{-1}\mathsf A^X_{1}(\theta)
\ldots\mathsf A^X_{n}(\theta)\right)}.
\notag \end{multline}

Consequently, up to multiplication by $\mathsf U^{g}_{\zzi}(\omega(x,\theta)^{-\frac c2})$
the matrices $\mathsf D^{X,g}_{\zzh}(\theta)$ can be obtained from elementary rotations.
The equalities above may also serve as alternative definitions.
In the cases $X=R,\frac R3,\frac L3,L$ one might even argue that
$\mathsf U^{g}_{\zzi}(\omega(x,\theta)^{-c})$ can be obtained from half-rotations.
\end{point}
\begin{point}
Every matrix $\mathsf M^g_{\zzi,\zzj}$ defined in this section can be twisted by rotations, so
we obtain the matrices $\mathsf M^g_{\zzi,\zzj}(x_1)=
\mathsf R_{\zzi}^{g,x_1}\mathsf M^g_{\zzi,\zzj}\mathsf R_{\zzj}^{g,-x_1}$.
For example, from $\mathsf M^g_{\zzi}(\theta)$  one can prepare the matrices
\[\mathsf M^g_{\zzi}(\theta,x_1)=\mathsf R_{\zzi}^{g,x_1}
\mathsf M^g_{\zzi}(\theta)\mathsf R_{\zzi}^{g,-x_1}.\]
In case of $\mathsf U_{\zzi}^g(x\mapsto f(x))$ we obtain $\mathsf U_{\zzi}^g(x\mapsto f(x-x_1))$ instead.
Then every arithmetical relation which holds among these kind of matrices remains true
for the extended matrices.
(But not the ones which contain transposition, because it is not invariant.)
So, in the computations above
``$(\theta)$'' should be replaced by ``$(\theta,x_1)$'', and
``$\omega(x,\theta)$''  should be replaced by ``$\omega(x-x_1,\theta)$''.
\end{point}
\begin{lemma}\label{lem:dupert} It yields
\begin{multline}
\mathsf F'^{X,g}_{\zzi}(\theta,x_1)
\left(\mathbf 1_{\zzi}\cos\tfrac{x_1-x_2}2+\mathsf H^{g}_{\zzi}\sin\tfrac{x_1-x_2}2  \right)
\mathsf F^{X,g}_{\zzi}(\theta,x_2)=\\
=\mathbf 1_{\zzi}\cos\tfrac{x_1-x_2}2(1+\sin^2\theta)
+\mathsf H^g_{\zzi}\sin\tfrac{x_1-x_2}2\cos^2\theta
-\mathsf U^g_{\zzi}\left(2\cos\left(\tfrac{x_1+x_2}2-x\right)\right)\sin\theta.
\notag\end{multline}
\qed\end{lemma}

\section{Deformations of multiplication actions}
\begin{point}\label{eq:expand}
Now we combine multiplication matrices and the special matrices introduced in the previous section.
Let us define
\[\mathsf U^{g}_{\zzi,\zzj}(a,\theta,x_1)=
\mathsf D^{g}_{\zzi}(\theta,x_1)
\mathsf U^{g}_{\zzi,\zzj}\left(a\right)
\mathsf D^{g}_{\zzj}(\theta,x_1)^{-1}.\]

Right now $\mathsf U^{g}_{\zzi,\zzj}(a,\theta,z_2)$ is defined only for
$\theta\notin \pi\zzh$, ie.~for $\sin\theta\neq\pm1$.
But we will shortly see that the definition extends.
Let us start with the case of $\mathsf U^g_{\zz}(a,\theta,x_1)$.
One finds that
\[\mathsf U^{g}_{\zz}\left(\sum_{n\in\zz}z^n a_n,\theta,x_1\right)=
\sum_{n\in\zz}\mathsf U^{g}_{\zz}(z^n,\theta,x_1)
\left([a_n]^+_g \mathbf 1_{\zz}+[a_n]^-_g\fli_{\zz}\right).\]

Using the abbreviations $s=\cos\theta$, $t=\sin\theta$, we see that for $n>0$ it yields
\[
\mathsf U^{g}_{\zz}(z^n,\theta,x_1)=
\left[\begin{array}{cccccc|c|cc}
\pmb\ddots&&&&&&&&\\
&1&&&&&&&\\\hline
&&s&tsz_1&\cdots&t^{n-1}sz_1^{n-1}&t^nz_1^n&&\\\hline
&&-tz_1^{-1}&s^2&\cdots&t^{n-2}s^2z_1^{n-2}&t^{n-1}sz_1^{n-1}&&\\
&&&\ddots&\ddots&\vdots&\vdots&&\\
&&&&\ddots&s^2&tsz_1&&\\
&&&&&-tz_1^{-1}&s&&\\
&&&&&&&1&\\
&&&&&&&&\pmb\ddots
\end{array}\right];
\]
and
\[\mathsf U_{\zz}^{g}(1,\theta,x_1)=\mathbf 1_{\zz};\]
furthermore
\[\mathsf U_{\zz}^{g}(z^{-n},\theta,x_1)=
\fli_{\zz}\mathsf U_{\zz}^{g}(z^n,\theta,-x_1)
\fli_{\zz},\]
which means that $\mathsf U^{g}_{\zz}(z^n,\theta,x_1)$
gets reflected and $z_1$ should be substituted by $z_1^{-1}$.

From the nature of these expressions we see that $\mathsf U^{g}_{\zz}(a,\theta,x_1 )$
extends to all $\theta$, and all along it yields a rapidly decreasing perturbation of
$\mathsf U^{g}_{\zz}(a)=\mathsf U^{g}_{\zz}(a,0,0)$.
The whole discussion extends to the general case.
In fact, due to multiplicativity,  it is sufficient to see that the matrices
\begin{align}
&\mathsf U^{g}_{\zz,\zzh}(z^{1/2},\theta,x_1 )=&
&\mathsf U^{g}_{\zzh,\zz}(z^{1/2},\theta,x_1 )=& \notag\\
&\quad=\left[
\begin{array}{cc|cc}
\pmb\ddots&&\\\hline
&\cos\frac\theta2&z_1\sin\frac\theta2&\\\hline
&-z_1^{-1}\sin\frac\theta2&\cos\frac\theta2&\\
&&&\pmb\ddots\\
\end{array}
\right],\quad& &\quad=\left[
\begin{array}{cc|c|c}
\pmb\ddots&&\\
&\cos\frac\theta2&z_1\sin\frac\theta2&\\\hline
&-z_1^{-1}\sin\frac\theta2&\cos\frac\theta2&\\
&&&\pmb\ddots\\
\end{array}
\right]& \notag
\end{align}
extend properly, which is clear.
Hence, from now on we define $\mathsf U^{g}_{\zzi,\zzj}(a,\theta,x_1)$
for all $\theta$ using the canonical continuous extension.
\end{point}
\begin{point}
We can similarly define
\[\mathsf U^{X,g}_{\zzi,\zzj}(a,\theta,x_1)=
\mathsf D^{X,g}_{\zzi}(\theta,x_1)
\mathsf U^{g}_{\zzi,\zzj}\left(a\right)
\mathsf D^{X,g}_{\zzj}(\theta,x_1)^{-1}.\]
Then the discussion is completely analogous.
For example, for $X=R$ it yields
\[
\mathsf U^{R,g}_{\zz}(z^n,\theta,x_1)=
\left[\begin{array}{cccccc|c|cc}
\pmb\ddots&&&&&&&&\\
&1&&&&&&&\\\hline
&&1&tz_1&\cdots&t^{n-1}z_1^{n-1}&t^nz_1^n&&\\\hline
&&-tz_1^{-1}&s^2&\cdots&t^{n-2}s^2z_1^{n-2}&t^{n-1}s^2z_1^{n-1}&&\\
&&&\ddots&\ddots&\vdots&\vdots&&\\
&&&&\ddots&s^2&ts^2z_1&&\\
&&&&&-tz_1^{-1}&s^2&&\\
&&&&&&&1&\\
&&&&&&&&\pmb\ddots
\end{array}\right];
\]
\begin{align}
&\mathsf U^{R,g}_{\zz,\zzh}(z^{1/2},\theta,x_1 )=&
&\mathsf U^{R,g}_{\zzh,\zz}(z^{1/2},\theta,x_1 )=&\notag\\
&=\quad\left[
\begin{array}{cc|cc}
\pmb\ddots&&&\\\hline
&1&&\\\hline
&-z_1^{-1}t&1&\\
&&&\pmb\ddots\\
\end{array}
\right],\quad&&\quad=\left[
\begin{array}{cc|c|c}
\pmb\ddots&&\\
&1&z_1t&\\\hline
&&1&\\
&&&\pmb\ddots\\
\end{array}
\right];&
\notag\end{align}
etc.
Note that $s^2=1-(tz_1)(tz_1^{-1})$.
This implies that in the cases $X=R,L$ the entries of matrices $\mathsf U_{\zzi,\zzj}^{X,g}(a,\theta,x_1)$
are (infinite) expressions of  $tz_1$, $tz^{-1}_1$ and the entries of $\mathsf U_{\zzi,\zzj}^g(a)$ only.
Similarly, in the cases $X=\frac R3,\frac L3$ the entries of matrices $\mathsf U_{\zzi,\zzj}^{X,g}(a,\theta,x_1)$
are  expressions of  $\tau z_1$, $\tau z^{-1}_1$ and the entries of $\mathsf U_{\zzi,\zzj}^g(a)$ only.
The various deformations are quite similar: Due to
\[\mathsf U^{Y,g}_{\zzi,\zzj}(a,\theta,x_1)=
\widetilde{\mathsf A}^{Y,X,g}_{\zzi}(\theta,x_1)
\mathsf U^{X,g}_{\zzi,\zzj}\left(a,\theta,x_1\right)
\widetilde{\mathsf A}^{X,Y,g}_{\zzj}(\theta,x_1)\]
(for $\theta\notin \pi\zzh$) and continuity we see that
the matrices corresponding to various choices of $X$ differ from
each other only in the environment of the coordinate axes.

\begin{lemma}
For $n\geq1$, $n\in\mathbb N$, it yields
\begin{multline}
\mathsf U^X_{\zz}(z^n,\theta,x_1)=
\mathsf A^X_{\frac12}(\theta,x_1)\mathsf A_{\frac32}^X(\theta,x_1)
\ldots\mathsf A^X_{n-\frac12}(\theta,x_1)\mathsf Z^n_{\zz}=\\=
\mathsf Z^n_{\zz}\mathsf A^X_{-n+\frac12}(\theta,x_1)
\ldots\mathsf A_{-\frac32}^X(\theta,x_1)\mathsf A^X_{-\frac12}(\theta,x_1),
\notag\end{multline}
\vspace{-6mm}
\begin{multline}
\mathsf U^X_{\zz,\zzh}(z^{n-\frac12},\theta,x_1)=
\mathsf A_{\frac12}^X(\theta,x_1) \ldots \mathsf A_{n-\frac32}^X(\theta,x_1) \daa^X_{n-\frac12}(\theta,x_1)
\mathsf Z^{n-\frac12}_{\zz,\zzh}=\\=
\mathsf Z^{n-\frac12}_{\zz,\zzh}\mathsf A_{-(n-1)}^X(\theta,x_1)\ldots
\mathsf A_{-1}^X(\theta,x_1) \daa^X_{0}(\theta,x_1),
\notag\end{multline}
\vspace{-6mm}
\begin{multline}
\mathsf U^X_{\zzh,\zz}(z^{n-\frac12},\theta,x_1)=
\daa^{-X}_{0}(-\theta,x_1)^{-1} \mathsf A_{1}^X(\theta,x_1)\ldots
\mathsf A_{n-1}^X(\theta,x_1)\mathsf Z^{n-\frac12}_{\zzh,\zz}=\\=
\mathsf Z^{n-\frac12}_{\zzh,\zz}
\daa^{-X}_{-n+\frac12}(-\theta,x_1)^{-1} \mathsf A_{-n+\frac32}^X(\theta,x_1)\ldots
\mathsf A_{-\frac12}^X(\theta,x_1),
\notag\end{multline}
\vspace{-6mm}
\begin{multline}
\mathsf U^X_{\zzh}(z^n,\theta,x_1)=
\daa^{-X}_{0}(-\theta,x_1)^{-1}\mathsf A_{1}^X(\theta,x_1)\cdot\ldots\cdot
\mathsf A_{n-1}^X(\theta,x_1)\daa^X_{n}(\theta,x_1)\mathsf Z^n_{\zzh}=\\=
\mathsf Z^n_{\zzh}\daa^{-X}_{-n}(-\theta,x_1)^{-1}
\mathsf A_{-n-1}^X(\theta,x_1)\cdot\ldots\cdot\mathsf A_{-1}^X(\theta,x_1)\daa^X_{0}(\theta,x_1).
\notag\end{multline}
\begin{proof}
It is sufficient to compute the case $n=1$;  the rest follows from the identities
$\mathsf Z^{n}\mathsf A_k^X(\theta,x_1)\mathsf Z^{-n} =\mathsf A_{k+n}^X(\theta,x_1)$,
$\mathsf Z^{n}\daa_k^X(\theta,x_1)\mathsf Z^{-n} =\daa_{k+n}^X(\theta,x_1)$.
\end{proof}
\end{lemma}
\end{point}
\begin{point}
It is quite obvious that the deformations $\mathsf U^{X,g}_{\zzi,\zzj}(a,\theta,x_1)$
can be constructed ``by hand'', without the help of the matrices $\mathsf D^g_{\zzi}(\theta,x_1)$.
In fact, this is sort of necessary.

In general, the deformations  $\mathsf U^{g}_{\zzi,\zzj}(a,\theta,x_1)$ are sort of the nicest
because they respect unitarity. Yet, the transformation parameter $\theta$ makes them look
sort of transcendental.
The deformations $\mathsf U^{R,g}_{\zzi,\zzj}(a,\theta,x_1)$ and
$\mathsf U^{L,g}_{\zzi,\zzj}(a,\theta,x_1)$ do not respect unitarity but they are of rather  algebraic
nature, hence general.
(They allow coefficients from more general rings.)
Their behaviour quite differs from the original case for $\theta=\frac\pi2$.
The deformations $\mathsf U^{\frac R3,g}_{\zzi,\zzj}(a,\theta,x_1)$ and
$\mathsf U^{\frac L3,g}_{\zzi,\zzj}(a,\theta,x_1)$ are sort of in the middle.
Still, they are probably closer to the plain case ($X=\emptis$).
\end{point}
\begin{point}
We see that the case $\theta=0$ means no deformation at all:
\[\mathsf U^{X,g}_{\zzi,\zzj}(a,0,x_1)=\mathsf U^g_{\zzi,\zzj}(a).\]

The other interesting value is $\theta=\frac\pi2$, where  very special kind of degeneracies occur:
\end{point}
\begin{lemma}\label{lem:ends}
It yields
\[\mathsf U^{g}_{\zzi,\zzj}(a,\tfrac{\pi}2,x_1)=
\mathsf L_{\zzi}^{g,\frac{x_1}2}
\left\{\begin{matrix}\mathsf U_{\zzi',\zzj'}^g(a) & \\& a(x_1) \end{matrix}\right\}^{g}_{\zzi,\zzj}
\mathsf L_{\zzj}^{g,-\frac{x_1}2},\]
\[\mathsf U^{\frac R3,g}_{\zzi,\zzj}(a,\tfrac{\pi}2,x_1)=
\mathsf L_{\zzi}^{g,\frac{x_1}2}
\left\{\begin{matrix}\mathsf U_{\zzi',\zzj'}^g(a) & \\& a(x_1) \end{matrix}\right\}^{R,g}_{\zzi,\zzj}
\mathsf L_{\zzj}^{g,-\frac{x_1}2},\]
\[\mathsf U^{\frac L3,g}_{\zzi,\zzj}(a,\tfrac{\pi}2,x_1)=
\mathsf L_{\zzi}^{g,\frac{x_1}2}
\left\{\begin{matrix}\mathsf U_{\zzi',\zzj'}^g(a) & \\& a(x_1) \end{matrix}\right\}^{L,g}_{\zzi,\zzj}
\mathsf L_{\zzj}^{g,-\frac{x_1}2},\]
\[\mathsf U^{R,g}_{\zzi,\zzj}(a,\tfrac{\pi}2,x_1)=
\mathsf L_{\zzi}^{g,\frac{x_1}2}
\left\{\begin{matrix}\mathsf U_{\zzi',\zzj'}^g(a) & \\ \bar{\mathsf S}(\Delta^{x_1}_{\zzj',\zzl}a)
& a(x_1) \end{matrix}\right\}^{R,g}_{\zzi,\zzj}
\mathsf L_{\zzj}^{g,-\frac{x_1}2},\]
\[\mathsf U^{L,g}_{\zzi,\zzj}(a,\tfrac{\pi}2,x_1)=
\mathsf L_{\zzi}^{g,\frac{x_1}2}
\left\{\begin{matrix}\mathsf U_{\zzi',\zzj'}^g(a) &{\mathsf S}(\Delta^{x_1}_{\zzi',\zzl}a)
\\& a(x_1) \end{matrix}\right\}^{L,g}_{\zzi,\zzj}
\mathsf L_{\zzj}^{g,-\frac{x_1}2}.\]
\qed\end{lemma}
\begin{point}For a multiplicatively invertible element $a\in\mathfrak B^l$ we define
\[\mathsf O^{L,g}_{\zzi}(a,\theta,x_1)=\mathbf 1_{\zzi}+t^2P^t_{\zz}\Bigg(x_1\mapsto
\mathsf L_{\zzi}^{g,\frac{x_1}2}
\left\{\begin{matrix}0&\mathsf S_{\zzi'}^{g}(\Xi^L_{\zzi',\zzl}a(x,x_1))\\&0\end{matrix}\right\}_{\zzi}^{L,g}
\mathsf L_{\zzi}^{g,-\frac{x_1}2}\Biggr),\]
\[\mathsf O^{R,g}_{\zzi}(a,\theta,x_1)=\mathbf 1_{\zzi}+t^2P^t_{\zz}\Bigg(x_1\mapsto
\mathsf L_{\zzi}^{g,\frac{x_1}2}
\left\{\begin{matrix}0&\\
\bar{\mathsf S}_{\zzi'}^{g}(\Xi^R_{\zzi',\zzl}a(x_1,x))&0\end{matrix}\right\}_{\zzi}^{R,g}
\mathsf L_{\zzi}^{g,-\frac{x_1}2}\Biggr),\]
where the Poisson transform acts in the variable $x_1$, entrywise in the matrices.

These are invertible matrices, and for $X=L,R$ we set
\[\mathsf U'^{X,g}_{\zzi,\zzj}(a,\theta,x_1)=\mathsf O^{X,g}_{\zzi}(a,\theta,x_1)^{-1}
\mathsf U^{X,g}_{\zzi,\zzj}(a,\theta,x_1).\]
Then

a.)
$\mathsf U'^{X,g}_{\zzi,\zzj}(a,\theta,x_1)$ is a smooth perturbation of
$\mathsf U'^{X,g}_{\zzi,\zzj}(a,0,x_1)=\mathsf U^{g}_{\zzi,\zzj}(a)$.

b.) In the special case $\theta=\frac\pi2$ it yields
\begin{multline}\mathsf U'^{L,g}_{\zzi,\zzj}(a,\tfrac{\pi}2,x_1)=
\mathsf L_{\zzi}^{g,\frac{x_1}2}
\left\{\begin{matrix}\mathbf 1_{\zzi'} &\mathsf S^g_{\zzi'}(-\frac 12\delta_{i,1} Ha_{\zzl}(x_1) a(x_1)^{-1} )\\& 1
\end{matrix}\right\}^{L,g}_{\zzi}\cdot \\ \cdot
\left\{\begin{matrix}\mathsf U_{\zzi',\zzj'}^g(a) &0
\\& a(x_1) \end{matrix}\right\}^{L,g}_{\zzi,\zzj}
\left\{\begin{matrix}\mathbf 1_{\zzj'} &\mathsf S^g_{\zzj'}(-\frac 12\delta_{j,1}
a(x)^{-1} H_{\zzl}a(x)
)\\& 1 \end{matrix}\right\}^{L,g}_{\zzj}
\mathsf L_{\zzj}^{g,-\frac{x_1}2}=
\\=\mathsf L_{\zzi}^{g,\frac{x_1}2}
\left\{\begin{matrix}\mathsf U_{\zzi',\zzj'}^g(a) & \mathsf S^g_{\zzi'}(
\widehat\Delta^{x_1}_{\zzi',\zzl}a(x) )
\\& a(x_1) \end{matrix}\right\}^{L,g}_{\zzi,\zzj}
\mathsf L_{\zzj}^{g,-\frac{x_1}2}.\notag\end{multline}
and
\begin{multline}\mathsf U'^{R,g}_{\zzi,\zzj}(a,\tfrac{\pi}2,x_1)=
\mathsf L_{\zzi}^{g,\frac{x_1}2}
\left\{\begin{matrix}\mathbf 1_{\zzi'} &\\ \bar{\mathsf S}^g_{\zzi'}
\left(\frac 12\delta_{i,1} a(x_1) H_{\zzl}(a^{-1})(x_1) \right)& 1
\end{matrix}\right\}^{R,g}_{\zzi}\cdot \\ \cdot
\left\{\begin{matrix}\mathsf U_{\zzi',\zzj'}^g(a) &0
\\& a(x_1) \end{matrix}\right\}^{R,g}_{\zzi,\zzj}
\left\{\begin{matrix}\mathbf 1_{\zzj'} &
\\\bar{\mathsf S}^g_{\zzj'}\left(\frac12\delta_{j,1}a(x) H_{\zzl}(a^{-1})(x)\right)& 1
\end{matrix}\right\}^{R,g}_{\zzj}
\mathsf L_{\zzj}^{g,-\frac{x_1}2}=\\
=\mathsf L_{\zzi}^{g,\frac{x_1}2}
\left(\left\{\begin{matrix}\mathsf U_{\zzj',\zzi'}^g(a^{-1}) &
\\\bar{\mathsf S}^g_{\zzi'}\left(\widehat\Delta^{x_1}_{\zzi',\zzl}(a^{-1})(x)\right)
& a(x_1)^{-1} \end{matrix}\right\}^{R,g}_{\zzj,\zzi}\right)^{-1}
\mathsf L_{\zzj}^{g,-\frac{x_1}2}.\notag\end{multline}

c.) The entries of  $\mathsf U'^{L,g}_{\zz}(a,\theta,x_1)$
are infinite algebraic expressions of  $tz_1$, $tz_1^{-1}$
and the entries of $\mathsf U^g_{\zzi,\zzj}(a)$ and  $\mathsf U^g_{\zzi,\zzj}(a)^{-1}$.
(This is not entirely trivial to check,
but it follows from Lemma \ref{lem:excorr} and Lemma \ref{lem:precorr}.)

d.) We remark that the correction terms
$\mathsf O^{X,g}_{\zzi}(a,\theta,x_1)$ are equal to $1_{\zzi}$
if $a$ is a constant or to an $\mathfrak A$-linear combination of $z^{1/2}$ and $z^{-1/2}$.
The matrices $\mathsf U'^{X,g}_{\zz}(a,\theta,x_1)$ are not multiplicative.
They are defined only to
yield a better values for $\theta=\frac\pi2$ than the ``unprimed'' versions.

In what follows we use the notation $\bix$ for upper left indices
when  the ``primed'' variations  $\prime L$ and  $\prime R$ are also allowed along the ordinary terms
$L,\frac L3,\emptis, \frac R3, R$.
\end{point}
\begin{point}\label{po:symmetries1}
In general, suppose that we have a matrix operation
which, given $a\in\mathfrak B^l$, $\theta$, $x_1$ it produces a $\zzi\times\zzj$ matrix
\[\mathsf M^{g}_{\zzi,\zzj}(a,\theta,x_1).\]
Such a construction can also be written down by the algebraic kernel functions
\[M_{\zzi,\zzj}(a,\theta,x_1;y_1,y_2).\]
It may or may not have certain symmetries.
Here are some of the possible symmetries described in the language of the matrices
and the equivalent transcription in terms of transformation kernels:

o.) Naturality in  $a$: If $\phi:\mathfrak A\rightarrow\mathfrak A'$ is a homomorphism then
\[\phi(\mathsf M^{g}_{\zzi,\zzj}(a,\theta,x_1))=\mathsf M^{g}_{\zzi,\zzj}(\phi(a),\theta,x_1);\]
\[\phi(M_{\zzi,\zzj}(a,\theta,x_1;y_1,y_2))=M_{\zzi,\zzj}(\phi(a),\theta,x_1;y_1,y_2).\]

a.) Rotation symmetry:
\[\mathsf M^{g}_{\zzi,\zzj}(R^{\beta} a,\theta,x_1)=
\mathsf R_{\zzi}^{g,\beta}\mathsf M^{X,g}_{\zzi,\zzj}(a,\theta,x_1-\beta)\mathsf R_{\zzj}^{g,-\beta};\]
\[M_{\zzi,\zzj}(R^\beta a,\theta,x_1;y_1,y_2)=M_{\zzi,\zzj}(a,\theta,x_1-\beta;y_1-\beta,y_2-\beta).\]

b.) Reflection  symmetry:
\[\mathsf M^{g}_{\zzi,\zzj}(Ca,\theta,x_1)=\fli_{\zzi}\mathsf M^{g}_{\zzi,\zzj}(a,\theta,-x_1)\fli_{\zzj};\]
\[M_{\zzi,\zzj}(Ca,\theta,x_1;y_1,y_2)=M_{\zzi,\zzj}(a,\theta,-x_1;-y_1,-y_2).\]

c.) Central symmetry:
\[\mathsf M^{g}_{\zzi,\zzj}(a,\theta,x_1)=\mathsf M^{g}_{\zzi,\zzj}(a,-\theta,x_1+\pi);\]
\[M_{\zzi,\zzj}(a,\theta,x_1;y_1,y_2)=M_{\zzi,\zzj}(a,-\theta,x_1+\pi;y_1,y_2).\]

d.) Central degeneracy:
\[\mathsf M^{g}_{\zzi,\zzj}(a,0,x_1-\beta)\text{ does not depend on $\beta$};\]
\[M_{\zzi,\zzj}(a,0,x_1-\beta;y_1,y_2)\text{ does not depend on $\beta$}.\]

(The last two conditions together mean that the associated matrix
depends essentially on $tz_1$, $tz_1^{-1}$ and $a$.)

e.) Duality symmetry:
If the associated matrices come in a variety depending on $X$ then a possible duality
symmetry  is
\[\mathsf M^{X,-g^\opp}_{\zzi,\zzj}(a^{\opp},\theta,x_1)=\mathsf M^{-X,g}_{\zzj,\zzi}(a,\theta,x_1)^{\top\opp}.\]
\[M^X_{\zzi,\zzj}(a^{\opp},\theta,x_1;y_1,y_2)=M^{-X}_{\zzj,\zzi}(a,\theta,x_1;y_2,y_1)^{\opp}.\]
(If there is no apparent choice for $X$ then it should be taken as the empty symbol.)

This discussion of symmetries applies even in the case when
$\mathsf M^{g}_{\zzi,\zzj}(a,\theta,x_1)$
would have less variables:
Indeed, the dependence of the construction  from $a,\theta,x_1$ may be nominal.
One can see that the constructions
$\mathsf U_{\zzi,\zzj}^g(a)$, $\mathsf F_{\zzi}^{X,g}(\theta,x_1)$, $\mathsf D_{\zzi}^{X,g}(\theta,x_1)$
and hence
$\mathsf U^{X,g}_{\zzi,\zzj}(a,\theta,x_1)$
have all those symmetries discussed above.
\end{point}

\begin{point}
For $a\in\mathfrak B^i$ we define
\[\mathsf E_{\zzi'}^g(a,x_1)=
\mathsf L^{g,\frac{x_1}2}_{\zzi'}
\left\{\begin{matrix}1_{\zzi}&\\&a(x_1)\end{matrix}\right\}_{\zzi'}^g,\qquad
\mathsf E_{\zzi'}^{\frac L3,g}(a,x_1)=
\mathsf L^{g,\frac{x_1}2}_{\zzi'}
\left\{\begin{matrix}1_{\zzi}&\\&a(x_1)\end{matrix}\right\}_{\zzi'}^{g,L},\]
\[\mathsf E_{\zzi'}^{L,g}(a,x_1)=
\mathsf L^{g,\frac{x_1}2}_{\zzi}
\left\{\begin{matrix}1_{\zzi}&\mathsf S_{\zzi}^g(\Delta^{x_1}_{\zzi,\zzi}a(x))
\\&a(x_1)\end{matrix}\right\}_{\zzi'}^{g,L},\]
\[\mathsf E'^{L,g}_{\zzi'}(a,x_1)=
\mathsf L^{g,\frac{x_1}2}_{\zzi'}
\left\{\begin{matrix}1_{\zzi}&\mathsf S_{\zzi}^g(-\frac12\delta_{i,0} H_{\zzi}a(x_1))
\\&a(x_1)\end{matrix}\right\}_{\zzi}^{g,L};\]
and
\[\widehat{\mathsf U}^g_{\zzi',\zzh}(a)=
\left\{\begin{matrix}\mathsf U^g_{\zzi,\zz}(a)&\\&1\end{matrix}\right\}_{\zzi',\zzh}^g,\qquad
\widehat{\mathsf U}^{\frac L3,g}_{\zzi',\zzh}(a)=\widehat{\mathsf U}^{L,g}_{\zzi',\zzh}(a)=
\left\{\begin{matrix}\mathsf U^g_{\zzi,\zz}(a)&\\&1\end{matrix}\right\}_{\zzi',\zzh}^{L,g},\]
\[\widehat{\mathsf U}'^{ L,g}_{\zzi',\zzh}(a)=
\left\{\begin{matrix}\mathsf U^g_{\zzi,\zz}(a)&\mathsf S_{\zzi}^g(-\frac12H_{\zzi}a(x))
\\&1\end{matrix}\right\}_{\zzi',\zzh}^{L,g}.\]
One can make analogous definitions with ``$R$''.
Then \[\mathsf U^{\bix,g}_{\zzi',\zzh}(a,\tfrac\pi2,x_1)=
\mathsf E^{\bix,g}_{\zzi'}(a,x_1)
\widehat{\mathsf U}^{\bix,g}_{\zzi',\zzh}(a)\mathsf L^{g,-\frac{x_1}2}_{\zzh}.\]

Furthermore, one can see that the matrices
\[\mathsf E^{\bix,g}_{\zzi'}(a,x_1)\mathsf L^{g,-\frac{x_1}2}_{\zzh}
\qquad\text{and}\qquad\mathsf L^{g,\frac{x_1}2}_{\zzh}
\widehat{\mathsf U}^{\bix,g}_{\zzi',\zzh}(a)\mathsf L^{g,-\frac{x_1}2}_{\zzh}\]
have naturality, and the rotation, reflection and duality symmetries.
\end{point}

\section{Kernel homotopies of rapidly decreasing matrices}
\begin{lemma}\label{lem:dupert0}
 For rapidly decreasing $\zzi\times \zzj$ matrices $\mathbf A^g$ the map
\[\mathsf P^{X,g}_{\zzi,\zzj}(\theta,x_1,x_2):\mathbf A^g\mapsto \mathsf F'^{X,g}_{\zzi}(\theta,x_2)^{-1}
\mathbf A^g \mathsf F^{X,g}_{\zzj}(\theta,x_1)^{-1}\cos^2\theta\]
extends smoothly to general $\theta$.
\begin{proof}
a.) Let us consider the case $X=\emptis$, $i=j=0$.
For the elementary matrix $\mathbf e_{n,m}$ ($n,m\geq0$) it yields
\[\mathbf e_{n,m}\mapsto
\left[
\begin{array}{cc|c|ccccc}
\ddots&&&&&\\
&0&&&&\\\hline
&&z_1^nz_2^{-m}1&z_1^nz_2^{-m+1}st&\ldots& z_1^nz_2^0st^m\\\hline
&&z_1^{n-1}z_2^{-m}st&z_1^{n-1}z_2^{-m+1}s^2t^2&\ldots& z_1^{n-1}z_2^0s^2t^{m+1}\\
&&\vdots&\vdots&&\vdots\\
&&z_1^1z_2^{-m}st^{n-1}&z_1^1z_2^{-m+1}s^2t^{n}&\ldots& z_1^1z_2^0s^2t^{n+m-1}\\
&&z_1^0z_2^{-m}st^n&z_1^0z_2^{-m+1}s^2t^{n+1}&\ldots& z_1^0z_2^0s^2t^{n+m}\\
&&&&&&0\\
&&&&&&&\ddots
\end{array}
\right]
\]
where $t=\sin\theta$, $s=\cos\theta$.
The behaviour is similar for $\mathbf e_{n,m}$ from other quadrants.
The result for  $a_{n,m}\mathbf e_{n,m}$ is slightly more complicated:
$a_{n,m}$ should be inserted between $z_1^p$ and $z_2^q$.
The other cases are all similar.
\end{proof}
\end{lemma}
\begin{point}
The case $\theta=0$ means no deformation:
 \[\mathsf P^{X,g}_{\zzi,\zzj}(0,x_1,x_2)(\mathbf A^g)=\mathbf A^g.\]
\end{point}

\begin{lemma}
If the rapidly decreasing matrix
$\mathbf A^g$ has  kernel function $A(x_1,x_2)$ then
\[\mathsf P^g_{\zzi,\zzj}(\tfrac\pi2,x_1,x_2)(\mathbf A^g)=
\mathsf L_{\zzi}^{g,\frac{x_1}2}
\left\{\begin{matrix}0& \\& A(x_1,x_2) \end{matrix}\right\}^{g}_{\zzi,\zzj}
\mathsf L_{\zzj}^{g,-\frac{x_2}2},\]
\[\mathsf P^{\frac R3,g}_{\zzi,\zzj}(\tfrac\pi2,x_1,x_2)(\mathbf A^g)=
\mathsf L_{\zzi}^{g,\frac{x_1}2}
\left\{\begin{matrix}0& \\& A(x_1,x_2) \end{matrix}\right\}^{R,g}_{\zzi,\zzj}
\mathsf L_{\zzj}^{g,-\frac{x_2}2},\]
\[\mathsf P^{\frac L3,g}_{\zzi,\zzj}(\tfrac\pi2,x_1,x_2)(\mathbf A^g)=
\mathsf L_{\zzi}^{g,\frac{x_1}2}
\left\{\begin{matrix}0& \\& A(x_1,x_2) \end{matrix}\right\}^{L,g}_{\zzi,\zzj}
\mathsf L_{\zzj}^{g,-\frac{x_2}2},\]
\[\mathsf P^{R,g}_{\zzi,\zzj}(\tfrac\pi2,x_1,x_2)(\mathbf A^g)=
\mathsf L_{\zzi}^{g,\frac{x_1}2}
\left\{\begin{matrix}0& \\\bar{\mathsf S}^g_{\zzj'}(\Delta^{x_2}_{\zzj',\zzj}A(x_1,x))&
A(x_1,x_2) \end{matrix}\right\}^{R,g}_{\zzi,\zzj}
\mathsf L_{\zzj}^{g,-\frac{x_2}2},\]
\[\mathsf P^{L,g}_{\zzi,\zzj}(\tfrac\pi2,x_1,x_2)(\mathbf A^g)=
\mathsf L_{\zzi}^{g,\frac{x_1}2}
\left\{\begin{matrix}0& \mathsf S^g_{\zzi'}(\Delta^{x_1}_{\zzi',\zzi}A(x,x_2))
\\& A(x_1,x_2) \end{matrix}\right\}^{L,g}_{\zzi,\zzj}
\mathsf L_{\zzj}^{g,-\frac{x_2}2}.\]
\qed\end{lemma}
\begin{lemma}
$\mathsf P^X_{\zz}(\theta,x_1,x_2)(u\mathbf e_{0,0})=u\mathbf e_{0,0}$.
\qed\end{lemma}

\section{Linearization of loops with base points}
\begin{point} For $a\in\mathfrak B^1$  we define
\[\mathsf K^{\bix,g}_{\zz}(a,\theta,x_1,x_2)=
\mathsf U^{\bix,g}_{\zz,\zzh}(a,\theta,x_1)\mathsf L_{\zzh}^{g,\frac {x_1-x_2}2}
\left(\mathsf U^{\bix,g}_{\zz,\zzh}(a,\theta,x_2)\right)^{-1},\]
and for $a\in\mathfrak B^0$  we define
\[\mathsf K^{\bix,g}_{\zzh}(a,\theta,x_1,x_2)=
\mathsf U^{\bix,g}_{\zzh}(a,\theta,x_1)\mathsf L_{\zzh}^{g,\frac {x_1-x_2}2}
\left(\mathsf U^{\bix,g}_{\zzh}(a,\theta,x_2)\right)^{-1}.\]
\end{point}
\begin{lemma} \label{lem:dupert1} It yields
\[\mathsf K^{X,g}_{\zzi}(a,\theta,x_1,x_2)=
\mathbf 1_{\zzi}\cos\tfrac {x_1-x_2}2+\left(\mathsf H_{\zzi}^g
+\mathsf P^X_{\zzi}(\theta,x_1,x_2)
\left(\mathsf J^g_{\zzi}(a)-\mathsf H_{\zzi}^g\right)\right)
\sin\tfrac {x_1-x_2}2.\]
\begin{proof}
First we prove the statement for  $\theta\notin\pi\zzh$.
From the definition of $\mathsf U^{X,g}_{\zzi,\zzj}(a,\theta,x_1)$, and
real multiplication commutativity with $b=\omega(x-x_1,\theta)$
we see that
\begin{multline}
\mathsf K^{X,g}_{\zzi'}(a,\theta,x_1,x_2)=
\mathsf F'^{X,g}_{\zzi'}(\theta,x_1)^{-1}\mathsf U^g_{\zzi',\zzh}(a)
\mathsf F'^{X,g}_{\zzh}(\theta,x_1)
\Bigl(\mathbf 1_{\zzh}\cos\tfrac {x_1-x_2}2+\\+\mathsf H_{\zzh}^g\sin\tfrac {x_1-x_2}2\Bigr)
\mathsf F^{X,g}_{\zzh}(\theta,x_2)\mathsf U^g_{\zzh,\zzi'}(a^{-1})
\mathsf F^{X,g}_{\zzi'}(\theta,x_2)^{-1}.
\notag\end{multline}

Applying  Lemma \ref{lem:dupert} in order expand the three terms in the middle,
applying real multiplication commutativity with  $b=\cos\left(\frac{x_1+x_2}2-x\right)$,
and applying  Lemma \ref{lem:dupert} again in order to write the terms back, we find
\begin{multline}\mathsf K^{X,g}_{\zzi'}(a,\theta,x_1,x_2)=
\mathbf 1_{\zzi'}\cos\tfrac {x_1-x_2}2+\mathsf H_{\zzi'}^g\sin\tfrac {x_1-x_2}2
+\\+\mathsf F'^{X,g}_{\zzi'}(\theta,x_1)^{-1}
\left(\mathsf J^g_{\zzi'}(a)-\mathsf H_{\zzi'}^g\right)
\mathsf F^{X,g}_{\zzi'}(\theta,x_1)^{-1}
\cos^2\theta\sin\tfrac {x_1-x_2}2.\notag\end{multline}
This proves the statement for  $\theta\notin\pi\zzh$.
Then, according to Lemma  \ref{lem:dupert0}, the result extends by continuity.
\end{proof}
\end{lemma}
\begin{theorem} \label{th:linear1}
 $\mathsf K^g_{\zz}$ satisfies the properties of a good linearizing cocycle:

i.) The cocycle property:
\[\mathsf K^{g}_{\zz}(a,\theta,x_1,x_1)=\mathbf 1_{\zz},
\qquad
\mathsf K^{g}_{\zz}(a,\theta,x_1,x_2)=
\mathsf K^{g}_{\zz}(a,\theta,x_1,x_3)\mathsf K^{g}_{\zz}(a,\theta,x_3,x_2).\]

ii.) $\mathsf K^{g}_{\zz}(a,\theta,x_1,x_2)$ is a smooth perturbation of
\[\mathbf 1_{\zz}\cos\tfrac {x_1-x_2}2+\mathsf H_{\zz}^{g}\sin\tfrac {x_1-x_2}2.\]

iii.) Values at $\theta=0$:
\[\mathsf K^{g}_{\zz}(a,0,x_1,x_2)=
\mathbf 1_{\zz}\cos\tfrac {x_1-x_2}2+\mathsf J_{\zz}^g(a)\sin\tfrac {x_1-x_2}2.\]

iv.) Values at $\theta=\tfrac\pi2$:
\[\mathsf K^{g}_{\zz}(a,\tfrac\pi2,x_1,x_2)=a(x_1)a(x_2)^{-1}\mathbf e_{0,0}
+\mathbf 1_{\zz\setminus \{0\}}\cos\tfrac {x_1-x_2}2+\mathsf H_{\zz}^{g}\sin\tfrac {x_1-x_2}2.\]

v.) If  $G$ is a skew-involution then
\[ \mathsf K^{g}_{\zz}(\cos\tfrac x2+G\sin\tfrac x2,\theta,x_1,x_2)=
\mathbf 1_{\zz}\cos\tfrac {x_1-x_2}2+\mathsf H_{\zz}^{g}[G]\sin\tfrac {x_1-x_2}2.\]
\begin{proof}
Points (i--iv) must be clear from the behaviour of $\mathsf U^g_{\zzi,\zzj}(a,\theta,x_1)$ discussed
previously.
Points (ii--v) follow from Lemma \ref{lem:dupert1} and the properties of the kernel homotopies.
\end{proof}
\end{theorem}
\begin{point} For
 $\mathsf K^{\bix,g}_{\zz}(a,\theta,x_1,x_2)$, in general,
the very same statement can be made except in point (iv) where in general
\[\mathsf K^{\bix,g}_{\zz}(a,\tfrac\pi2,x_1,x_2)=
\mathsf E^{\bix,g}_{\zz}(a,x_1)\mathsf E^{\bix,g}_{\zz}(a,x_2)^{-1}.\]
Actually, this makes difference only for
\[\mathsf K^{L,g}_{\zz}(a,\tfrac\pi2,x_1,x_2)=\left\{
\begin{matrix}\mathsf L^{g,\frac{x_1-x_2}2}_{\zzh}&
\mathsf L^{g,\frac{x_1}2}_{\zzh}\mathsf S^g_{\zzh}(\Delta^{x_1}_{\zzh,\zzh}a(x)
-\Delta^{x_2}_{\zzh,\zzh}a(x))a(x_2)^{-1}\\&a(x_1)a(x_2)^{-1}\end{matrix}
\right\}_{\zz}^{X,g},\]
and for $\mathsf K^{R,g}_{\zz}(a,\tfrac\pi2,x_1,x_2)$, where a similar formula holds.
Nevertheless, these are only very mild deformations of $\mathsf K^{g}_{\zz}(a,\tfrac\pi2,x_1,x_2)$;
the difference is only a mostly inconsequential off-diagonal term.

There is an analogous statement for $\mathsf K^{\bix,g}_{\zzh}(a,\theta,x_1,x_2)$, in general,
but it is  less nice in (iv), trivial in (v), and its intuitive meaning is not clear.
In particular, it does not a provide a naive linearization procedure for periodic loops.
\end{point}
\begin{point}
It is worthwhile to examine the invariance properties of the linearizing cocycles constructed above.
Compared to the discussion in \ref{po:symmetries1} the difference is that we have two
parameter variables: $x_1$ and $x_2$.
The generalizations of naturality, the rotation, reflection, central symmetries
and central degeneracy
are straightforward, it is simply required to double the variable $x_1$:
Along every expression $f(x_1)$ an other one, $f(x_2)$, must be inserted.
However,

e'.) Duality symmetry is better to be formulated, in general, as
\[\mathsf M^{X,-g^\opp}_{\zzi,\zzj}(a^{\opp},\theta,x_1,x_2)
=\mathsf M^{-X,g}_{\zzj,\zzi}(a^{-1},\theta,x_2,x_1)^{\top\opp},\]
\[M_{\zzi,\zzj}(a^{\opp},\theta,x_1,x_2;y_1,y_2)=M_{\zzj,\zzi}(a^{-1},\theta,x_2,x_1;y_2,y_1)^{\opp}.\]
in terms of matrices and kernels respectively.

Then  $\mathsf L_{\zzh}^{g,\frac{x_1-x_2}2}$, and, more generally,
$\mathsf K_{\zzi}^{X,g}(a,\theta,x_1,x_2)$ have all these symmetries.
\end{point}
\begin{point}
Even at this point it becomes apparent that it might be useful to formulate the linearization
procedure entirely in terms of cocycles.
For a loop $a$ we might form the trivial cocycle
\[(x,\tilde x)\mapsto \tilde a(x,\tilde x)=a(x)a(\tilde x)^{-1}.\]
If $(x,\tilde x)\rightarrow c(x,\tilde x)$  is a cocycle then we define its dual cocycle by
\[(x,\tilde x)\mapsto c^{\dual}(x,\tilde x)=c(\tilde x,x)^\opp.\]

e''.) If we have a matrix construction $\mathsf M^{X,g}_{\zzj,\zzi}(x,\theta,x_1,x_1)$
depending on a cocycle $c$ then we define duality invariance by
\[\mathsf M^{X,-g^\opp}_{\zzi,\zzj}(c^{\dual},\theta,x_1,x_2)
=\mathsf M^{-X,g}_{\zzj,\zzi}(c,\theta,x_2,x_1)^{\top\opp}.\]
\[M_{\zzi,\zzj}(c^{\dual},\theta,x_1,x_2;y_1,y_2)=M_{\zzj,\zzi}(c,\theta,x_2,x_1;y_2,y_1)^{\opp}.\]
in terms of matrices and kernels respectively.

First, it is clear from the real representation that $\mathsf K_{\zzi}^{X,g}(a,\theta,x_1,x_2)$
depends only on the cocycle $\tilde a$ and not on its trivialization $a$.
Hence, we could have written $\mathsf K_{\zzi}^{X,g}(\tilde a,\theta,x_1,x_2)$.
Then we may observe that it satisfies duality invariance in the above sense.
In fact, we can define rotated cocycles by $R^\beta c(x_1,x_2)=c(x_1-\beta,x_2-\beta)$, etc.,
and we could have formulated rotation invariance, etc., in these terms.
\end{point}

If we consider linearization of loops in strict sense then this corresponds to finding a
trivialization of the linearizing cocycle.
If we choose such a trivialization then we generally loose some kind of
invariance property.
\begin{point}
When we consider loops as pointed maps, ie. when $a(0)=1$ holds then
\[\widetilde{\mathsf K}_{\zz}^{g}(a,\theta,x_1)=\mathsf K_{\zz}^{g}(a,\theta,x_1,0)\]
yields a linearizing homotopy which leaves that class of loops invariant.
It yields a trivialization of the cocycle in the sense that
\[\mathsf K_{\zz}^{g}(a,\theta,x_1,x_2)=
\widetilde{\mathsf K}_{\zz}^{g}(a,\theta,x_1)\widetilde{\mathsf K}_{\zz}^{g}(a,\theta,x_2)^{-1} .\]

Here both the trivial loop extension
\[\widetilde{\mathsf K}^{g}_{\zz}(a,\tfrac\pi2,x_1)=
a(x_1)\mathbf e_{0,0}+\mathbf 1_{\zz\setminus \{0\}}\cos\tfrac {x_1}2+\mathsf H_{\zz}^{g}\sin\tfrac {x_1}2\]
and the linearized loop
\[\widetilde{\mathsf K}^{g}_{\zz}(a,0,x_1)
=\mathbf 1_{\zz}\cos\tfrac {x_1}2+\mathsf J_{\zz}^g(a)\sin\tfrac {x_1}2\]
are nice.
Its nice properties are immediate from the statement of Theorem \ref{th:linear1}.
However, there is no apparent rotation invariance in this construction.
There are natural variants $\widetilde{\mathsf K}_{\zz}^{\bix,g}(a,\theta,x_1)$ of this construction, of course.
\end{point}

\begin{point}
\label{po:NN}
We have another construction which is less invariant than $\mathsf K^g_{\zz}(a,\theta,x_1,x_2)$
but in many ways better than  $\widetilde{\mathsf K}^g_{\zz}(a,\theta,x_1)$:

Let us consider
\[\mathsf N_{\zz}^{g}(a,\theta,x_1)=
\mathsf U^{g}_{\zz,\zzh}(a,\theta,x_1)
\mathsf L_{\zzh}^{g,\frac{x_1}2}
\widehat{\mathsf U}^{g}_{\zz,\zzh}(a)^{-1}
=\mathsf U^{g}_{\zz,\zzh}(a,\theta,x_1)
\mathsf U^{g}_{\zz,\zzh}(a,\tfrac\pi2,x_1)^{-1}
\mathsf E_{\zz}^{g}(a,x_1).\]
It has the trivializing property
\[\mathsf K_{\zz}^{g}(a,\theta,x_1,x_2)=
\mathsf N_{\zz}^{g}(a,\theta,x_1)\mathsf N_{\zz}^{g}(a,\theta,x_2)^{-1} .\]

The advantage of this function is that it directly linearizes
\[\mathsf N_{\zz}^{g}(a,\tfrac\pi2,x_1)=a(x_1)\mathbf e_{0,0}
+\mathbf 1_{\zz\setminus \{0\}}\cos\tfrac {x_1}2+\mathsf H_{\zz}^{g}\sin\tfrac {x_1}2,\]
although the linearized loop
\[\mathsf N_{\zz}^{g}(a,0,x_1)=\mathsf U^{g}_{\zz,\zzh}(a)
\left(\mathbf 1_{\zzh}\cos\tfrac {x_1}2+\mathsf H_{\zzh}^g\sin\tfrac {x_1}2\right)
\widehat{\mathsf U}^{g}_{\zz,\zzh}(a)^{-1}
\] is not particularly nice in itself.
Nevertheless, the failing of the construction to be rotation invariant
is not worse than that of $\mathsf L^{g,\frac{x_1}2}_{\zzh}$ itself:
the equation of rotation invariance holds up to multiplication by  $\mathsf L^{g,\frac{\beta}2}_{\zzh}$.
What we can say that here linearization ``acts'' on the left side while the right side is
undetermined up to multiplication by a constant.
One can consider variants $\mathsf N_{\zz}^{\bix,g}(a,\theta,x_1)$
which may or may not have better properties (like algebraic parametrizability in terms of $t$, etc.).

We remark that in the special cases $a=z^{\pm1/2}$ it yields
\[\mathsf N^{X,g}_{\zz}(z^{1/2},\theta,x_1)
=\mathsf L^{g,x_1}_{\zz}[z^{1/2}]\daa^X_{\frac12}(\theta)\daa^X_{\frac12}(\tfrac\pi2)^{-1},  \]
\[\mathsf N^{X,g}_{\zz}(z^{-1/2},\theta,x_1)
=\mathsf L^{g,x_1}_{\zz}[z^{-1/2}]\daa^{-X}_{-\frac12}(-\theta)\daa^{-X}_{-\frac12}(-\tfrac\pi2)^{-1}.  \]
(The primed versions do not differ in these cases.)
\end{point}
\begin{point}
Of the various constructions above the cases $\bix=R,L,\prime R,\prime L$
allow easy algebraization, ie.~more general coefficient rings.
In these cases $\mathsf U^{\bix,g}_{\zzi,\zzj}(a,\theta,x_1)$,
$\mathsf E^{\bix,g}_{\zz}(a,x_1)$ and $\widehat{\mathsf U}^{\bix,g}_{\zz,\zzh}(a)$
are expressed in terms of $tz_1$, $tz_1^{-1}$, and the entries of $\mathsf U_{\zzi}^g(a)$;
not even division by $2$ is required.
\end{point}

\section{Non-commutative cyclic loops}
Having gained some experience with linearizing loops we can try to generalize the results above.
The main technique we use is simply suppression of variables, which allows to formulate statements
in greater generality.

\begin{point}
Consider a superalgebra $\mathfrak A=\mathfrak A^+\oplus\mathfrak A^-$, and consider an extension
$\mathfrak A[\mathsf z,\mathsf z^{-1}]$
where that $\mathsf z$ is a cyclic variable  which does not necessarily
commute with elements of $\mathfrak A$ but
$\mathsf z\mathfrak A^+\mathsf z^{-1}\subset\mathfrak A^+,$
$\mathsf z^{-1}\mathfrak A^+\mathsf z\subset \mathfrak A^+$ and
$\mathsf z\mathfrak A^-\mathsf z\subset\mathfrak A^-,$
$\mathsf z^{-1}\mathfrak A^-\mathsf z^{-1}\subset \mathfrak A^-$.
Then $\mathfrak A[\mathsf z,\mathsf z^{-1}]$
is naturally $\bo Z_2\times \mathbb Z$ graded, and assume that it is endowed
by a (not necessary complete) locally convex algebra topology and
$\mathfrak B^{0,0}$ is the algebra of rapidly decreasing sequences regarding the grading.
Ie. what  we consider smooth loops with a  non-commuting loop variable.
One wonders if some kind of linearization is possible in this case.
\end{point}
\begin{point}
We can introduce the symbol $\mathsf z^{1/2}$, and we can take
$\mathfrak B^0=\mathfrak B^{0,0}\oplus \mathfrak B^{1,1}=\mathfrak A[\mathsf z,\mathsf z^{-1}]\oplus
\mathsf z^{1/2}\mathfrak A[\mathsf z,\mathsf z^{-1}]\mathsf z^{1/2}$
and
$\mathfrak B^1=\mathfrak B^{0,1}\oplus \mathfrak B^{1,0}=
\mathfrak A[\mathsf z,\mathsf z^{-1}]\mathsf z^{1/2}\oplus
\mathsf z^{1/2}\mathfrak A[\mathsf z,\mathsf z^{-1}]$.
Regarding linearizations, what we have to do is to provide analogous formulas with a suppression of variables.
That means we have to replace ``$g,x_1$'' by $\mathsf z$ and
``$g,\frac{x_1}2$'' by $\mathsf z^{1/2}$ as much as it possible.
But we have to keep in mind that, say, $\mathfrak A^+$ is not invariant for conjugation
by $\mathsf z^{1/2}$.

If $a\in\mathfrak B^{i,j}$ then according to the choice of basis $\mathbf e_n\leftrightarrow\mathsf z^n$
the representing matrix of  multiplication by $a$ on the left is a matrix
$\mathsf U^{\mathsf z}_{\zzi,\zzj}(a)$.
Let  $a=\sum_{n\in\zzi}\mathsf z^na_n$ where $a_n\in\mathfrak A$ or
$\mathsf  z^{1/2}\mathfrak A^+\mathsf  z^{-1/2}\oplus \mathsf  z^{1/2}\mathfrak A^-\mathsf  z^{1/2}$.
The effect of multiplication by $a$ on the left is given as
$a \mathsf z^k= \sum_{n\in\zzi} \mathsf z^{n+k} (\mathsf z^{-k}a_n^+\mathsf z^k)+
\mathsf z^{n-k} (\mathsf z^{k}a_n^-\mathsf z^k)$.
According to the correspondence of bases $\mathbf e_n\leftrightarrow\mathsf z^n$
the representing matrix of  multiplication by $a$ on the left is
\[\mathsf U^{\mathsf z}_{\zzi,\zzj}(a)=
\sum_{n\in\zzi,k\in\zzj} \left(\mathsf z^{-k}a_{n-k}^+\mathsf z^k
+\mathsf z^{k}a_{n+k}^-\mathsf z^k\right)  \mathbf e_{n,k}.\]
Symbolic rotation actions are given by
$\mathsf R^{\mathsf z}_{\zzi}=\sum_{n\in\zzi}\mathsf z^{-n}\mathbf e_{n,n}.$

Now, we say that a $\zzi\times\zzj$ matrix $A$ is of Toeplitz type if
cutting it up to four pieces by a vertical and a horizontal line we get Toeplitz matrices.
We say that $A$ is of rotationally Toeplitz type if
$(\mathsf R^{\mathsf z}_{\zzi})^{-1}A(\mathsf R^{\mathsf z}_{\zzj})$
is of Toeplitz type.
We can notice that
$(\mathsf R^{\mathsf z}_{\zzi})^{-1}
\mathsf U^{\mathsf z}_{\zzi,\zzj}(a)
(\mathsf R^{\mathsf z}_{\zzj})
= \sum_{n\in\zzi,k\in\zzj} \mathsf z^{n-k}a_{n-k}\,  \mathbf e_{n,k} $.
So, the multiplication actions are at least of rotationally Toeplitz type.
That makes topology clear: We will use the Toeplitz topology induced by a symbolic rotation,
and that is compatible to the original loop topology.
\end{point}
\begin{point}
We can proceed further by other definitions like
$\mathsf L_{\zz}^{\mathsf z}[u]=\mathsf z\mathbf 1_{\zz^-}+u\mathbf e_{0,0}+\mathsf z^{-1}\mathbf 1_{\zz^+}$,
$\mathsf L_{\zzh}^{\mathsf z}=\mathsf z\mathbf 1_{\zzh^-}
+\mathsf z^{-1}\mathbf 1_{\zzh^+}$, etc.
Proceeding analogously to the skew-complex case then we see that
$\mathsf U^{X,\mathsf z}_{\zzi,\zzj}(a,\theta)$
will also be of rotationally Toeplitz type, etc.
It may be strange that the expression $\mathsf U^{X,\mathsf z}_{\zzi,\zzj}(a,\theta)$ would correspond to
$\mathsf U^{X,g}_{\zzi,\zzj}(a,\theta,x)$, ie.~$x$  appears as a variable of $a(x)$ and
also appears as an external parameter (as we have an apparent shortage of variables).
This, in fact, causes no problems, the first occurrence of $x$ is closed relative to the second.
It is similar to expressions of type
$\exists x\mathcal F(x,x_1,\ldots, x_n;\exists x\mathcal G(x,x_1,\ldots,x_n))$ in logic.

Regarding the special case $\theta=\frac\pi2$ we cannot write expressions which are similarly
compact (as we cannot use $g$, we cannot simply take difference functions, etc.), nevertheless the structure of
the matrices remains highly special,  analogous to the skew-complex case.
Ultimately, taking direct analogues of the various matrices as $\zzi\times\zzj$ matrices
(and not as block decompositions) we can define $\widehat{\mathsf U}^{\bix,\mathsf z}_{\zzi',\zzh}(a)$ etc.
In theory this requires that we expand all our matrices as in, say, point \ref{eq:expand} and make sure
that we can proceed.
For $a\in\mathfrak B^{0,1}$ we find that $\mathsf U_{\zz,\zzh}^{X,\mathsf z}(a,\theta)$ is matrix
with coefficients  from $\mathfrak B^{0,1}$ and
$\widehat{\mathsf U}^{\bix,\mathsf z}_{\zzh}(a)$ is a matrix with coefficients from
$\mathsf z^{1/2}\mathfrak A^+\mathsf z^{-1/2}\oplus \mathsf z^{1/2}\mathfrak A^-\mathsf z^{1/2}$.
Then we can indeed present a linearization:
\end{point}
\begin{point}
For $ a\in\mathfrak B^{0,1}$ we can define
\[\mathsf N_{\zz}^{\mathsf z}(a,\theta)=\mathsf U^{\mathsf z}_{\zz,\zzh}(a,\theta)
\mathsf L_{\zzh}^{\mathsf z^{1/2}}\widehat{\mathsf U}^{g}_{\zz,\zzh}(a)^{-1}
=\mathsf U^{\mathsf z}_{\zz,\zzh}(a,\theta)
\mathsf U^{\mathsf z}_{\zz,\zzh}(a,\tfrac\pi2)^{-1}
\mathsf L_{\zz}^{\mathsf z^{1/2}}[a(\mathsf z)].\]
\end{point}
\begin{theorem} It yields:

i.) $\mathsf N_{\zz}^{\mathsf z}(a,\theta)$ is rapidly decreasing perturbation of
$\mathsf L_{\zz}^{\mathsf z^{1/2}}[0]$ with entries from $\mathfrak B^{0,1}$.

ii.)
$\mathsf N_{\zz}^{\mathsf z}(a,0)=\mathsf U^{\mathsf z}_{\zz,\zzh}(a)
\mathsf L_{\zzh}^{\mathsf z^{1/2}}\widehat{\mathsf U}^{\mathsf z}_{\zzh,\zz}(a)^{-1}$
has entries from $\mathfrak A\mathsf z^{1/2}\oplus\mathfrak A\mathsf z^{-1/2}$.

iii.)
$\mathsf N_{\zz}^{\mathsf z}(a,\tfrac\pi2)=\mathsf L^{\mathsf z^{1/2}}_{\zz}[a(\mathsf z)].$

iv.) In the special cases $a=\mathsf z^{\pm1/2}$ it gives
$\mathsf N_{\zz}^{\mathsf z}(\mathsf z^{1/2},\theta)=
\mathsf L^{\mathsf z^{1/2}}[\mathsf  z^{1/2}]\daa_{\frac12}(\theta)\daa_{\frac12}(\tfrac\pi2)^{-1}  ,$
and
$\mathsf N_{\zz}^{\mathsf z}(\mathsf z^{-1/2},\theta)=
\mathsf L^{\mathsf z^{1/2}}[\mathsf  z^{-1/2}]
\daa_{\frac12}(-\theta)\daa_{\frac12}(-\tfrac\pi2)^{-1}$.
\qed\end{theorem}
Other variants of this linearization procedure exist naturally.
We can linearize  loops $a\in\mathfrak  B^{0,0}$ in more conventional sense by taking, for example,
$\mathsf N_{\zz}^{\mathsf z}(a\mathsf z^{-1/2},\theta)\mathsf z^{1/2}$.

\section{Linearization in terms of transformation kernels}
\begin{point}
Let us consider loops and matrix actions on them using Fourier expansions with respect to $g=\mathrm i$.
Let
$\mathbf a=\sum_{k\in\zzj} a_k \mathbf e_k$,
$\mathbf C=\sum_{n\in\zzi,m\in\zzj} c_{n,m}\mathbf e_{n,m}$
be rapidly decreasing matrices.
Through the correspondences
\[\mathbf a=\sum_{k\in\zzj} a_k \mathbf e_k\leftrightarrow a(x)=\sum_{k\in\zzj} a_k\mathrm e^{\mathrm i kx}\]
\[\mathbf C=\sum_{n\in\zzi,m\in\zzj} c_{n,m}\mathbf e_{n,m}\leftrightarrow C(x,\tilde x)
=\sum_{n\in\zzi,m\in\zzj} c_{n,m} \mathrm e^{\mathrm i (nx-m\tilde x)} \]
we see that various matrix operations can be represented according to
\[\mathbf C\mathbf a \leftrightarrow C*a(x)=\int_{y=0}^{2\pi} C(x,y)a(y)\frac{dy}{2\pi},\]
\[\mathbf C_1\mathbf C_2\leftrightarrow C_1*C_2(x,\tilde x)=
\int_{y=0}^{2\pi} C_1(x,y)C_2(y,\tilde x) \frac{dy}{2\pi}, \]
\[\mathbf C^\top \leftrightarrow C^{\mathrm T}(x,\tilde x)=C(-\tilde x,- x).\]
As  transformation kernels do not depend on any choice like $g=\mathrm i$
we see that the operations ``$*$'' and ``$\mathrm T$'' are entirely canonical.
After this integration will be understood over $[0,2\pi]$
(or over a compactly supported smooth distribution such that the sum of its translates by
$2\pi\mathbb Z$ is the constant $1$ function.)

Using these observations, in what follows, we will describe our linearizations
in terms of transformation kernels.
Instead of writing a transformation kernel as $C(x_1,x_2)$ we will sometimes write it as
$C[x_1,\alpha]$ where $\alpha=x_1-x_2$. Then
\[ C[x_1,\alpha]=C(x_1,x_1-\alpha) \qquad=\qquad C(x_1,x_2)=C[x_1,x_1-x_2].\]

In this formalism
\[C*a(x)=\int_\beta C[x,\beta] a(x-\beta)\frac{d\beta}{2\pi}
=\int_\beta C[x,\beta] R^\beta a(x)\frac{d\beta}{2\pi},\]
\[ C_1*C_2[x,\alpha]=\int_\beta C_1[x,\beta]C_2[x-\beta,\alpha-\beta]\frac{d\beta}{2\pi}
=\int_\beta C_1[x,\beta]R^\beta C_2[x,\alpha-\beta]\frac{d\beta}{2\pi},\]
\[C^{\mathrm T}[x,\alpha]= C[\alpha-x,\alpha];\]
here $R^\alpha$ means rotated in the first variable.
\end{point}

\begin{point}
Actually, the most natural operations cannot be represented by smooth transformation kernels.
We can allow transformation kernels which are distributions in the variable $\alpha$,
and there is a smooth dependence in $x$.
For $b\in\mathfrak B^j$ the standard multiplication and rotation actions are given by
\[U_{\zzi,\zzj}(a)*b(x) =a(x)b(x)= \int_\alpha a(x)\delta_{\zzj}^{0}(\alpha)\, R^\alpha b(x),\]
\[R_{\zzj}^\beta*b(x)= \int_\alpha \delta_{\zzj}^{\beta}(\alpha)\, R^\alpha b(x);\]
some Dirac $\delta$-functions appear.
Hence we find that the kernels are
\[U_{\zzi,\zzj}(a)[x,\alpha]= a(x)\delta_{\zzj}^{0}(\alpha)\quad=\quad
U_{\zzi,\zzj}(a)(x,\tilde x)= a(x)\delta_{\zzj}^{0}(x-\tilde x)= a(\tilde x)\delta_{\zzi}^{0}(x-\tilde x),\]
\[R_{\zzj}^\beta[x,\alpha]=\delta_{\zzj}^{\beta}(\alpha)\quad=\quad
R_{\zzj}^\beta[x,\tilde x]=\delta_{\zzj}^{\beta}(x-\tilde x).\]

For a transformation kernel $C$ representing a $\zzi\times\zzj$ matrix
one finds that the following identities are computationally useful:
\[C*U_{\zzj,\zz_k}(a)[x,\alpha]=C[x,\alpha] a(x-\alpha)\quad=\quad
C*U_{\zzj,\zz_k}(a)(x,\tilde x)=C(x,\tilde x) a(\tilde x),\]
\[U_{\zz_h,\zzi}(a)*C[x,\alpha]=a(x)C[x,\alpha]\quad=\quad
U_{\zz_h,\zzi}(a)*C(x,\tilde x)=a(x)C(x,\tilde x),\]
\[R_{\zzi}^\beta*C*R_{\zzj}^{-\beta}[x,\alpha]=C[x-\beta,\alpha]\quad=\quad
R_{\zzi}^\beta*C*R_{\zzj}^{-\beta}(x,\tilde x)=C(x-\beta,\tilde x-\beta).\]

In particular, a transformation is rotation invariant if and only if
its kernel  depends only on $\alpha$.
\end{point}
\begin{point}
Such transformation are, for example, the Poisson, Poisson-Hilbert, and Hilbert transformations.
Using the notation $w=\mathrm e^{\mathrm i\alpha}$ we find that their kernels are
\[P^r_{\zz}[x,\alpha]=\frac{1-r^2}{(1-rw)(1-rw^{-1})}=\frac{1-r^2}{1-2r\cos \alpha+r^2},\]
\[H^r_{\zz}[x,\alpha]=\frac{-\mathrm ir(w-w^{-1})}{(1-rw)(1-rw^{-1})}
=\frac{2r\sin\alpha}{1-2r\cos \alpha+r^2}.\]
\[H_{\zz}[x,\alpha]=\PV\mathrm i\frac{w+1}{w-1}=\PV \cot\frac\alpha2.\]
\[P^r_{\zzh}[x,\alpha]=\frac{r^{1/2}(1-r)(w^{1/2}+w^{-1/2})}{(1-rw)(1-rw^{-1})}
=\frac{2r^{1/2}(1-r)\cos\frac\alpha2}{1-2r\cos \alpha+r^2},\]
\[H^r_{\zzh}[x,\alpha]=\frac{-\mathrm ir^{1/2}(1+r)(w^{1/2}-w^{-1/2})}{(1-rw)(1-rw^{-1})}
=\frac{2r^{1/2}(1+r)\sin\frac\alpha2}{1-2r\cos \alpha+r^2}\]
\[H_{\zzh}[x,\alpha]=\PV\frac{2\mathrm i}{w^{1/2}-w^{-1/2}}=\PV \cosec\frac\alpha2.\]

So, we can  forget the first variable $x$.
After this we omit the sign $\PV$.
\end{point}
\begin{point}
The Dirac delta and Cauchy principal values are among the very mildest singularities.
Yet, heuristical computations with them should be taken with a grain of salt.
For example, when one takes the square of the  Hilbert transform one finds
\[\cot\frac\beta2\cot\frac{\alpha-\beta}2-1=
\begin{cases}
-\cosec^2\frac\beta2&\text{if }\alpha\in2\pi\zz,\\
\cot\frac\alpha2\left(\cot\frac\beta2+\cot\frac{\alpha-\beta}2\right)&\text{if }\alpha\notin2\pi\zz,\\
\end{cases}\]
and
\[\cosec\frac\beta2\cosec\frac{\alpha-\beta}2=
\begin{cases}
-\cosec^2\frac\beta2&\text{if }\alpha\in4\pi\zz,\\
{\cosec\frac\alpha2}\left(\cot\frac\beta2+\cot\frac{\alpha-\beta}2\right)&\text{if }\alpha\notin2\pi\zz,\\
{\cosec^2\frac\beta2}&\text{if }\alpha\in4\pi\zzh,\\
\end{cases}\]
from which
\[\int\cot\frac\beta2\cot\frac{\alpha-\beta}2 \,\frac{d\beta}{2\pi} -1=
-\delta_{\zz}^0(\alpha)\]
and
\[\int\cosec\frac\beta2\cosec\frac{\alpha-\beta}2 \, \frac{d\beta}{2\pi}=
-\delta_{\zzh}^0(\alpha). \]
are quite believable but it needs some faith regarding the singular values.

If one wants to be more precise, it is better to rely on the following lemma
in order to approximate properly:
\begin{lemma}\label{lem:hilbex} For $0\leq r_2<r_1$
\begin{multline}
P^{r_1}_{\zzi}(\beta)P^{r_2}_{\zzi}(\alpha-\beta)=\\=
\frac{P^{r_1}_{\zz}(\beta)+P^{r_2}_{\zz}(\alpha-\beta)}2\,P^{r_1r_2}_{\zz_i}(\alpha)
+\frac{P^{r_1}_{\zz}(\beta)-P^{r_2}_{\zz}(\alpha-\beta)}2\,P^{r_2/r_1}_{\zz_i}(\alpha)
\\
-\frac{H^{r_1}_{\zz}(\beta)+H^{r_2}_{\zz}(\alpha-\beta)}2\,H^{r_1r_2}_{\zz_i}(\alpha)
+\frac{H^{r_1}_{\zz}(\beta)+H^{r_2}_{\zz}(\alpha-\beta)}2\,H^{r_2/r_1}_{\zz_i}(\alpha),
\notag\end{multline}
\begin{multline}
P^{r_1}_{\zzi}(\beta)H^{r_2}_{\zzi}(\alpha-\beta)
=\\=
\frac{P^{r_1}_{\zz}(\beta)+P^{r_2}_{\zz}(\alpha-\beta)}2\,H^{r_1r_2}_{\zz_i}(\alpha)
+\frac{P^{r_1}_{\zz}(\beta)-P^{r_2}_{\zz}(\alpha-\beta)}2\,H^{r_2/r_1}_{\zz_i}(\alpha)
\\
+\frac{H^{r_1}_{\zz}(\beta)+H^{r_2}_{\zz}(\alpha-\beta)}2\,P^{r_1r_2}_{\zz_i}(\alpha)
-\frac{H^{r_1}_{\zz}(\beta)+H^{r_2}_{\zz}(\alpha-\beta)}2\,P^{r_2/r_1}_{\zz_i}(\alpha),
\notag\end{multline}
\begin{multline}
H^{r_1}_{\zzi}(\beta)P^{r_2}_{\zzi}(\alpha-\beta)
=\\=
\frac{P^{r_1}_{\zz}(\beta)+P^{r_2}_{\zz}(\alpha-\beta)}2\,H^{r_1r_2}_{\zz_i}(\alpha)
+\frac{-P^{r_1}_{\zz}(\beta)+P^{r_2}_{\zz}(\alpha-\beta)}2\,H^{r_2/r_1}_{\zz_i}(\alpha)
\\
+\frac{H^{r_1}_{\zz}(\beta)+H^{r_2}_{\zz}(\alpha-\beta)}2\,P^{r_1r_2}_{\zz_i}(\alpha)
+\frac{H^{r_1}_{\zz}(\beta)+H^{r_2}_{\zz}(\alpha-\beta)}2\,P^{r_2/r_1}_{\zz_i}(\alpha),
\notag\end{multline}
\begin{multline}H^{r_1}_{\zzi}(\beta)H^{r_2}_{\zzi}(\alpha-\beta)-\delta_{i,0}\cdot 1
=\\=
-\frac{P^{r_1}_{\zz}(\beta)+P^{r_2}_{\zz}(\alpha-\beta)}2\,P^{r_1r_2}_{\zz_i}(\alpha)
+\frac{P^{r_1}_{\zz}(\beta)-P^{r_2}_{\zz}(\alpha-\beta)}2\,P^{r_2/r_1}_{\zz_i}(\alpha)
\\
+\frac{H^{r_1}_{\zz}(\beta)+H^{r_2}_{\zz}(\alpha-\beta)}2\,H^{r_1r_2}_{\zz_i}(\alpha)
+\frac{H^{r_1}_{\zz}(\beta)+H^{r_2}_{\zz}(\alpha-\beta)}2\,H^{r_2/r_1}_{\zz_i}(\alpha),
\notag\end{multline}
and similarly with $\beta$ and $\alpha-\beta$ interchanged. \qed
\end{lemma}
\end{point}
\begin{point}
If a transformation kernel also depends on an additional parameter $x_1$ then
we can write it in alternative ways, as
\[C(x_1,y_1,y_2)=C[x_1,y_1,\beta]=C[[x_1,\eta,\beta]],\]
where $\beta=y_1-y_2$, $\eta=x_1-y_1$.
In this case, for the multiplication it yields
\begin{multline}
C_1*C_2(x_1,y_1,y_2)=\int_{u}C_1(x_1,y_1,u)C_2(x_1,u,y_2)\dfrac{du}{2\pi}=
\\=C_1*C_2[x_1,y_1,\beta]=\int_{\chi}C_1[x_1,y_1,\chi]C_2[x_1,y_1-\chi,\beta-\chi]\dfrac{d\chi}{2\pi}=
\\=C_1*C_2[[x_1,\eta,\beta]]
=\int_{\chi}C_1[[x_1,\eta,\chi]]C_2[[x_1,\eta+\chi,\beta-\chi]]\dfrac{d\chi}{2\pi}.
\notag\end{multline}
\end{point}

\begin{point}
From the kernel of the Hilbert transform we find that for $a\in\mathfrak B^i$
\[J_{\zzi'}(a)(y_1,y_2)=\frac{a(y_1)a(y_2)^{-1}}{\sin\frac{y_1-y_2}2}=
\frac{\tilde a(y_1,y_2)}{\sin\frac{y_1-y_2}2}.\]

Let us concentrate to the case $X=\emptis$. We find
\[F_{\zz}(\theta,x_1)=\delta_{\zz}^0(y_1-y_2)(1-\cos ( y_1-x_1)\sin\theta)
+\frac{\sin(\frac{y_1+y_2}2-x_1)}{\sin\frac{y_1-y_2}2}\sin\theta+(\cos\theta-1),\]
\[F'_{\zz}(\theta,x_2)=\delta_{\zz}^0(y_1-y_2)(1-\cos (y_2-x_2)\sin\theta)
-\frac{\sin( \frac{y_1+y_2}2-x_2 )}{\sin\frac{y_1-y_2}2}\sin\theta+(\cos\theta-1).\]

From Lemma \ref{lem:dupert1} we know that
\begin{multline}
\mathsf K^{g}_{\zz}(a,\theta,x_1,x_2)=\cos\frac{x_1-x_2}2\mathbf 1_{\zz}+
\sin\frac{x_1-x_2}2\biggl(\mathsf H^g_{\zz}+\\+
\mathsf F^g_{\zz}(\theta,x_1)\mathsf U^g_{\zz}(\omega(x-x_1,\theta)^{-1})
(\mathsf J^g_{\zz}(a)-\mathsf H^g_{\zz})\mathsf U^g_{\zz}(\omega(x-x_2,\theta)^{-1})
\mathsf F'^g_{\zz}(\theta,x_2)\cos^2\theta\biggr).
\notag\end{multline}
\end{point}
\begin{point}\label{po:speckernel}
Then we can write down the transformation kernel
\begin{multline}K_{\zz}(a,\theta, x_1,x_2,y_1,y_2)
=\cos\frac{x_1-x_2}2\delta^0_{\zz}(y_1-y_2)
+\sin\frac{x_1-x_2}2\Biggl(\cot\frac{y_1-y_2}2+\\
+\Biggl[\frac{1-\cos (y_1-x_1)\sin\theta}{\omega(y_1-x_1,\theta)}
\Biggl(\frac{\tilde a(y_1,y_2)-\cos\frac{y_1-y_2}2}{\sin\frac{y_1-y_2}2}\Biggr)
\frac{1-\cos (y_2-x_2)\sin\theta)}{\omega(y_2-x_2,\theta)}+\\
+\int_{u_1}\left(\frac{\sin(\frac{y_1+u_1}2-x_1)}{\sin\frac{y_1-u_1}2}\sin\theta+(\cos\theta-1)
\right)\Biggl(\frac{\tilde a(u_1,y_2)-\cos\frac{u_1-y_2}2}{
\omega(u_1-x_1,\theta)\sin\frac{u_1-y_2}2}\Biggr)\frac{du_1}{2\pi}
\cdot\\\cdot\frac{1-\cos (y_2-x_2)\sin\theta}{\omega(y_2-x_2,\theta)}+\ldots
\Biggr]\cos^2\theta\Biggr)\notag.
\end{multline}
It has the following properties:

i.) It is a smooth perturbation of
\[\cos\frac{x_1-x_2}2 \delta_{\zz}^0\left(y_1-y_2\right)
+\sin\frac{x_1-x_2}2 \cot\frac{y_1-y_2}2.\]

ii.) It has the cocycle property:
\[K(a,\theta,x_1,x_1,y_1,y_2)=\delta_{\zz}^0 (y_1-y_2),\]
\[K_{\zz}(a,\theta,x_1,x_2,y_1,y_2)=
\int_{u}  K(a,\theta,x_1,x_3,y_1,u)K(a,\theta,x_3,x_2,u,y_2). \]

iii.) Value at $\theta=0$:
\[K(a,0,x_1,x_2,y_1,y_2)=
\cos\tfrac{x_1-x_2}2 \delta_{\zz}^0(y_1-y_2)
+\sin\tfrac{x_1-x_2}2\cosec\tfrac{y_1-y_2}2\,\tilde a(y_1,y_2).\]

iv.) Value at $\theta=\frac\pi2$:
\[K(a,\tfrac\pi2,x_1,x_2,y_1,y_2)=
\cos\tfrac{x_1-x_2}2 \left(\delta_{\zz}^0(y_1-y_2)-1\right)+
\sin\tfrac{x_1-x_2}2\cot\tfrac{y_1-y_2}2+\tilde a(x_1,x_2).\]

v.) If $G$ is a skew-involution then
\[K(\cos\tfrac\lambda2+G\sin\tfrac\lambda2,\theta,x_1,x_2,y_1,y_2)=
\cos\tfrac{x_1-x_2}2\delta_{\zz}^0(y_1-y_2)+
\sin\tfrac{x_1-x_2}2\left(\cot\tfrac{y_1-y_2}2+G\right).\]

The various symmetries are manifested in the kernel.
For example, rotation symmetry is related to the fact that
the transformation can be expressed relatively simply by the cocycle $\tilde a$, $x_1$, and
$\alpha=x_1-x_2$, $\eta=x_1-y_1$, $\beta=y_1-y_2$.
\end{point}

\section{Linearization of  $\mathbb S^1$-cocycles}
\begin{point}
\textit{Fourier expansion of  $\mathbb S^1$-actions.}
Assume that $\mathfrak B=\mathfrak B^0\oplus \mathfrak B^1$ is a $\bo Z_2$-graded algebra.
Assume that it allows a smooth one-parameter group of automorphisms
$R:\mathbb R\times \mathfrak B\rightarrow \mathfrak B$ which is
periodic on $\mathfrak B^0$ and skew-periodic on   $\mathfrak B^1$.
We use $R^\alpha$ to denote the action of $R$ belonging to the parameter value $\alpha$.
Then the periodicity conditions are $R^{2\pi}b=b$ for $b\in\mathfrak B^0$ and
$R^{2\pi}b=-b$ for $b\in\mathfrak B^1$.

This situation is characterized by the lack of an explicit loop variable $x$.
Nevertheless, even in this case one can apply Fourier decomposition in terms of the rotation action.
Indeed, due to smoothness the action of $R$ allows a Fourier decomposition
\[R^\alpha a =\frac12 a_{[0]}+\sum_{n\in\zz^+} a_{[n]}\cos n\alpha+ H_{\zz}a_{[n]}\sin n\alpha
\qquad\text{for}\qquad a\in\mathfrak B_0,\]
\[R^\alpha a =\sum_{n\in\zzh^+} a_{[n]}\cos n\alpha+ H_{\zzh}a_{[n]}\sin n\alpha
\qquad\text{for}\qquad a\in\mathfrak B_1;\]
where the Hilbert transforms are defined by
\[H_{\zz}a=\int_{\beta} \cot\frac\beta2 R^\beta a \,\frac{d\beta}{2\pi}\qquad\text{and}\qquad
H_{\zzh}a=\int_{\beta} \cosec\frac\beta2 R^\beta a \,\frac{d\beta}{2\pi};\]
and more particularly,
\[a_{[n]}=2\int_{\beta} \cos n\beta\, R^\beta a \,\frac{d\beta}{2\pi}\qquad\text{and}\qquad
H_{\zz+n}a_{[n]}=2\int_{\beta} \sin n\beta\, R^\beta a \,\frac{d\beta}{2\pi}.\]

(Cf.: If $b\in\mathfrak B^0=\mathcal C_{\zz}(\mathfrak A)$ and
\[b(x)=b_0+\sum_{n\in\zz^+}\left(b_n\cos nx+b_{-n}\sin nx\right)\]
then
\[R^\alpha b(x)= b_0+\sum_{n\in\zz^+}(b_n\cos n(x-\alpha)+b_{-n}\sin n(x-\alpha))=\]
\[= b_0+\sum_{n\in\zz^+}\left(b_n\cos nx+b_{-n}\sin nx\right)\cos n\alpha
+\left(b_n\sin nx-b_{-n}\cos nx\right)\sin n\alpha= \]
\[= b_0+\sum_{n\in\zz^+}\left(b_n\cos nx+b_{-n}\sin nx\right)\cos n\alpha
+H_{\zz}\left(y\mapsto b_n\cos ny+b_{-n}\sin ny\right)(x)\sin n\alpha. \]
Similar formula holds for $b\in\mathfrak B^1$.)

Applying the auxiliary involution $\mathrm i$ we can take complex Fourier expansion, too.
This expansion can be used as follows.
\end{point}

\begin{point}
\textit{The spaces $\mathcal K_{i,j}(\mathfrak B,R)$ and $\Psi_{i,j}(\mathfrak B,R)$.}
Let $\mathcal K_{i,j}(\mathfrak B,R)$ denote the space $\mathcal C_{\zzj}(\mathfrak B^l)$ of
smooth  $\mathfrak B^l$-valued loops on $\mathbb R$, which are periodic for $j=0$ and skew-periodic for $j=1$
(in what follows: $[j]$-periodic) endowed with  product operations as follows:
For $b\in\mathfrak B^j$ and $C\in\mathcal K_{i,j}(\mathfrak B,R)$ we define
\[C*b=\int_\beta C(\beta) R^\beta b\, \frac{d\beta}{2\pi}\in\mathfrak B^i.\]
So, $C$ represents a map $\mathfrak B^j\rightarrow \mathfrak B^i$.
For $C_1\in \mathcal K_{i,j}(\mathfrak B,R)$  and $C_2\in \mathcal K_{j,k}(\mathfrak B,R)$ we set
\[C_1*C_2(\alpha)=\int_\beta C_1(\beta) R^\beta C_2(\alpha-\beta)\, \frac{d\beta}{2\pi}.\]
A natural rotation action is given by
\[(R^\gamma C)(\alpha)= R^\gamma(C(\alpha)).\]

In this situation $C$ allows a Fourier expansion
\[R^\beta C(\alpha)=\sum_{n\in\zzi,m\in\zzj}
c_{n,m}\mathrm e^{\mathrm i(-n\beta+m(\beta+\alpha))},\qquad c_{n,m}\in\mathfrak B^l.\]
Then $C*$ can be represented by the matrix $\sum c_{n,m}\mathbf e_{n,m}$.
This example shows that the usual matrix computations can be transferred to these cases
where ``the variables are not explicit''.

If, instead of smooth loops, we use $\mathfrak B^l$-valued distributions then
one can see that the multiplication operation above still exist
(it is only a matter of partial integration).
In that way we obtain the much larger algebra $\Psi_{i,j}(\mathfrak B,R)$.
What we really need is $\widetilde{\mathcal K}_{i,j}(\mathfrak B,R)$
which is the extension of $\mathcal K_{i,j}(\mathfrak B,R)$ by the
identity and Hilbert transform elements $1_{\zzi}=\delta^0_{\zzi}(\alpha)$ and
$H_{\zzi}=h^0_{\zzi}(\alpha)$.
In case of the Hilbert transforms it is particularly convenient to
use matrices after we have established their properties using Lemma \ref{lem:hilbex}
or in some other ways.
\end{point}
\begin{point}
\textit{Cocycles.}
We say that $C\in\mathcal C_i(\mathfrak B^0)$ is a cocycle if the identities
\[C(0)=1\qquad\text{and}\qquad C(\alpha)=C(\beta)R^\beta C(\alpha-\beta)\]
hold. We denote the set of such cocycles by $\CC_i(\mathfrak B^0)$.

To any invertible element $a\in\mathfrak B^j$ we can associate the trivial cocycle
\[\tilde a(\alpha)=a R^\alpha a^{-1},\]
which is in $\CC_i(\mathfrak B^0)$.
Not all cocycles are trivializable. For example, there may be cocycles in
$\CC_1(\mathfrak B^0)$ even if $\mathfrak B^1=0$.
In fact, as skew-periodic cocycles are $\mathfrak B^0$-valued loops we see that
skew-periodic cocycles have, a priori,  nothing  to do with skew-periodic loops.
A second remark is that in the case of $\mathbb C$ coefficients, in general,
one cannot construct a skew-periodic loop from a periodic one naturally, but
one can construct a skew-periodic cocycle from a periodic one naturally.

A skew-periodic cocycle is linear, if it is of shape
\[C(\beta)=\cos\frac\beta2+G\sin\frac\beta2,\]
where $G$ is a skew-involution and $G$ is $R$-invariant (ie. $R^\alpha G=G$ for all $\alpha$).
\end{point}

\begin{point}\textit{``Matrix'' extension.}
Let $\mathcal K_{i,j}(\mathfrak B^k)$
contain those $\mathfrak B^k$-valued smooth functions
\[(\eta,\beta) \mapsto C(\eta,\beta)\]
which are $[i+j]$-periodic in $\eta$ and $[j]$-periodic in $\beta$, and which multiply as
\[C_1*C_2(\eta,\beta)=\int_\gamma C_1(\eta,\gamma)C_2(\eta+\gamma,\beta-\gamma)\frac{d\gamma}{2\pi}.\]
Then $\Psi_{i,j}(\mathfrak B^k)$  can be defined similarly,
these are distributions in the variable $\beta$ and smooth in the variable $\eta$; etc.
Passing to the variable $-\eta$ instead of $\eta$ we see that these algebras are indeed
isomorphic to matrix algebras.

It may not be obvious but these spaces also behave as tensor products:
One can construct maps
\[\mathfrak B^{i+j+k}\otimes \Psi_{i,j}(\mathfrak B,R)\xrightarrow{e}
\Psi_{i,j}(\mathfrak B^k)\]
where
\[e(b \otimes A)(\eta,\beta) =b R^{\eta}(A(\beta)).\]
One can see that via this construction both $\mathfrak B^s$ and $\Psi_{i,j}(\mathfrak B)$ embed
through taking tensor products with $1$.
(We remark, however, that $\mathcal K_{1,0}(\mathfrak B^0)\neq 0$ is possible
even if $\mathfrak B^1\neq0$. Hence the image of the map above not necessarily dense.)

There are two kinds of rotations to be considered here:
One is the external rotation
\[(R^\mu C)(\eta,\beta)= R^\mu(C(\eta,\beta)), \]
and the other one is the internal rotation
\[(\widehat R^\mu C)(\eta,\beta) = R^\mu(C(\eta-\mu,\beta)).\]

They are related to the tensor product construction by the commutative diagrams
\[\xymatrix{ \mathfrak B^{i+j+k}\otimes  \Psi_{i,j}(\mathfrak B,R)\ar[r]\ar[d]_{R^\gamma\otimes R^\gamma}&
\Psi_{i,j}(\mathfrak B^k)\ar[d]^{R^\gamma}\\
\mathfrak B^{i+j+k}\otimes  \Psi_{i,j}(\mathfrak B,R)\ar[r]&
\Psi_{i,j}(\mathfrak B^k)}
\quad ~_{\begin{matrix}~\\\text{and}\end{matrix}}\quad
\xymatrix{ \mathfrak B^{i+j+k}\otimes\Psi_{i,j}(\mathfrak B,R)\ar[r]\ar[d]_{R^\gamma\otimes \Id}&
\Psi_{i,j}(\mathfrak B^k)\ar[d]^{\widehat R^\gamma}\\
\mathfrak B^{i+j+k}\otimes\Psi_{i,j}(\mathfrak B,R)\ar[r]&
\Psi_{i,j}(\mathfrak B^k).}\]

When we are interested in finding cocycles in that algebra, then the cocycle variable
will be $\alpha$ and a candidate cocycle
will be given by the function $C(\alpha,\eta,\beta)$;
dependence on additional factors will be indicted in variables in front of $\alpha$.
Then the $\Psi_{i,i}(\mathfrak B^0)$-valued linear cocycles
associated to the internal rotations are of shape
\[\cos\frac\alpha2+\sin\frac\alpha2 R^\eta G(\beta),\]
where $G(\beta)\in \Psi_{i,i}(\mathfrak B,R)$ is an involution.
\end{point}
We can formulate our general theorem:
It is related to the special kernel in \ref{po:speckernel} by
taking the cocycle $c(x_1,\lambda)=\tilde a(x_1,x_1-\lambda)$,
$\alpha=x_1-x_2$, $\eta=x_1-y_1$, $\beta=y_1-y_2$ and suppressing the variable $x_1$.
\begin{theorem} We can define $K_{\zz}(c,\theta,\alpha,\eta,\beta)$ such that

i.) It is a smooth perturbation of
\[\cos\frac{\alpha}2 \delta_{\zz}^0(\beta)+\sin\frac{\alpha}2 h^0_{\zz}(\beta).\]

ii.) It has the cocycle property with respect to internal rotations:
\[K(c,\theta,0,\eta,\beta)=\delta_{\zz}^0(\beta),\]
\[K_{\zz}(c,\theta,\alpha,\eta,\beta)=
\int_{\gamma}  K(c,\theta,\sigma,\eta,\gamma)
R^{\sigma}K(c,\theta,\alpha-\sigma,\eta-\sigma+\gamma,\beta-\gamma)\frac{d\gamma}{2\pi}. \]

iii.) Value at $\theta=0$:
\[K(c,0,\alpha,\eta,\beta)=
\cos\frac\alpha2 \delta_{\zz}^0(\beta)+\sin\frac{\alpha}2h^0_{\zzh}(\beta)\,R^{\eta}c(\beta),\]

iv.) Value at $\theta=\frac\pi2$:
\[K(a,\tfrac\pi2,\alpha,\eta,\beta)=
\cos\frac\alpha2 \left(\delta_{\zz}^0(\beta)-1\right)+\sin\frac{\alpha}2h^0_{\zz}(\beta)+c(\alpha).\]

v.) If $G$ is a skew-involution
\[K(\cos\tfrac\lambda2+G\sin\tfrac\lambda2,\theta,\alpha,\eta,\beta)=
\cos\frac\alpha2\delta_{\zz}^0(\beta)+
\sin\frac{\alpha}2\left(h^0_{\zz}(\beta)+G\right).\]

vi.) It is equivariant with respect to external rotations:
\[ K_{\zz}(R^\gamma a,\theta,\alpha,\eta,\beta)=
R^\gamma K_{\zz}(a,\theta,\alpha,\eta,\beta).\]
\begin{proof}
It is analogous to the matrix proof of loop linearization case but
using the quadruple Fourier expansion in $R^\xi,\alpha,\eta,\beta$.
\end{proof}
\end{theorem}
\begin{point}
In the abstract language by that we have defined a map
\[K_{\zz}:\CC_1(\mathfrak B^0)\times\mathbb R\rightarrow
\CC_1(\widetilde{\mathcal K}_{0,0}(\mathfrak B^0))\]
with nice properties, it linearizes (from $\theta=\frac\pi2$ to $\theta=0$)
the trivial extension of the cocycle $c$
to the linear cocycle associated to the Bott skew-involution
\[J_{\zz}(c)=\left(\beta\mapsto h^0_{\zzh}(\beta) c(\beta)\right)
\in\widetilde{\mathcal K}_{0,0}(\mathfrak B,R). \]
\end{point}

\begin{point}
Let us take a look to concrete example where the above theorem have meaning:
Suppose that $P\rightarrow U$ is an smooth, oriented $\mathbb S^2$-bundle
where the spheres come with a standard metric.
Assume that $V$ is a smooth $K(P)$-representative which means for us that a sufficiently nice
 splitting of $\mathbb C_{\zz}=\bigoplus_{k\in\zz}\mathbb C\mathbf e_{k}$
is given over every point of $P$.
Assume that $s:X\rightarrow P$ is a section, and $V|_{s(X)}$ is trivial.
Let $W$ be the vertical, $\mathbb S^2$-tangent bundle of $P$ over $s(X)$.

Then on the unit circle $\mathbb S^1W_{s(u)}$ of a fiber $W_{s(u)}$  we can canonically
construct a $\GL^\infty(\mathbb C_{\zz})$-valued clutching cocycle for $V$
 through exponentiation and parallel transport.
Due to the orientedness of the $\mathbb S^2$-bundle there is a canonical
$\mathbb S^1$-action on $\mathbb S^1W$, hence it
yields a periodic $\mathbb S^1$-cocyle.
Due to the presence of complex numbers we can canonically construct a
skew-periodic cocycle.
Then, through linearization, this is homotopic to a nice
$\GL^\infty(\mathbb C_{\zzh}\otimes \mathcal C (\mathbb S^1W_{s(u)}))$-valued cocycle.
This is determined by a sufficiently nice skew-involution, or rather splitting of
in  $\mathbb C_{\zz}\otimes \mathcal C_0 (\mathbb S^1W_{s(u)})$.
Due to the presence of many space-dimensions we see that these are homotopic to
splittings where the nontivial part of the splitting happens in $\mathbb C_{\zz}\otimes
\mathbb C 1_{s(u)}.$
We can realize that this yields the  Thom isomorphism
$K(P,X)\rightarrow K(X)$ in this limited setting.
In fact, this was a very special situation, we applied the linearization theorem
to principal $\mathbb S^1$-bundles instead of general $\mathbb S^1$-actions, etc.
\end{point}

\begin{point}
If one wants to linearize trivial cocycles (trivially) then
one can define $N(a,\theta,\eta,\beta)$ similarly.
\end{point}

\end{document}